\documentclass{amsart}

\usepackage{url}
\usepackage{tikz}
\usepackage{eso-pic}
\usepackage{listings}
\usepackage{booktabs}
\usepackage{subfig}

\theoremstyle{plain}
\newtheorem{theorem}{Theorem}
\newtheorem*{theorem*}{Theorem}
\newtheorem{lemma}[theorem]{Lemma}
\newtheorem*{lemma*}{Lemma}

\newtheorem{corollary}[theorem]{Corollary}

\theoremstyle{definition}
\newtheorem{definition}{Definition}
\newtheorem*{definition*}{Definition}
\newtheorem{observation}[theorem]{Observation}

\newtheorem{example}{Example}

\newcommand{\RR}{\ensuremath{\mathbb R}}

\newcommand{\NN}{\ensuremath{\mathbb{N}}}

\newcommand{\G}{\ensuremath{\mathcal{G}}}
\newcommand{\g}{\ensuremath{\mathfrak{g}}}
\newcommand{\0}{\ensuremath{\vec 0}}
\newcommand{\x}{\ensuremath{\bar x}}
\newcommand{\y}{\ensuremath{\bar y}}
\newcommand{\Fld}{\ensuremath{\mathbb{F}}}
\newcommand{\class}[1]{\ensuremath{\overline{#1}}}
\newcommand{\bmat}{\begin{bmatrix}}
\newcommand{\emat}{\end{bmatrix}}
\newcommand{\st}{\,\mid\,}

\DeclareMathOperator{\rank}{rank}

\DeclareMathOperator{\mr}{mr}
\DeclareMathOperator{\characteristic}{char}
\DeclareMathOperator{\diag}{diag}
\DeclareMathOperator{\projdim}{pdim}

\def \fh{\mathrm{fullhouse}}

\begin{document}

\author{Jason Grout}
\address{Department of Mathematics, Iowa State University, Ames, Iowa 50011}

\email{grout@iastate.edu}


\title{The minimum rank problem over finite fields}

\begin{abstract}
  The structure of all graphs having minimum rank at most $k$ over a
  finite field with $q$ elements is characterized for any possible $k$
  and $q$.  A strong connection between this characterization and
  polarities of projective geometries is explained.  Using this
  connection, a few results in the minimum rank problem are derived by
  applying some known results from projective geometry.
\end{abstract}

\maketitle

\section{Introduction}

Given a field $F$ and a simple undirected graph $G$ on $n$ vertices (i.e., an
undirected graph without loops or multiple edges), let
$S(F,G)$ be the set of symmetric $n\times n$ matrices $A$ with entries
in $F$ satisfying $a_{ij} \neq 0$, $i \neq j$, if and only if $ij$ is
an edge in $G$.  There is no restriction on the diagonal entries of
the matrices in $S(F,G)$.  Let
\begin{equation*}
\mr(F,G)=\min\{\rank A \st A \in S(F,G)\}.
\end{equation*}
Let $\G_k(F)=\{G \st \mr(F,G)\leq k\}$, the set of simple graphs with
minimum rank at most $k$.  

The problem of finding $\mr(F,G)$ and describing $\G_k(F)$ has
recently attracted considerable attention, particularly for the case
in which $F = \RR$ (see
\cite{nylen-minrank,cdv-eig-mult,johnson-duarte-trees, hsieh-minrank,
  johnson-saiago-maxmult, chen, vdh-nullity, BFH1-minrankpath,
  barrett-vdHL-minrank2-infinite, hall, arav,
  bento-duarte-tridiag-matrices, bfh-mult-pathcover-tree, bfh-cdv,
  barrett-vdHL-minrank2-finite, ding-kotlov-minrank-finite,
  barioli-fallat-joins}).  The minimum rank problem over $\RR$ is a
sub-problem of a much more general problem, the inverse eigenvalue
problem for symmetric matrices: given a family of real numbers, find
every symmetric matrix that has the family as its eigenvalues.  More
particularly, the minimum rank problem is a sub-problem of the inverse
eigenvalue problem for graphs, which fixes a zero/nonzero pattern for
the symmetric matrices considered in the inverse eigenvalue problem.
The minimum rank problem can also be thought of in this way: given a
fixed pattern of off-diagonal zeros, what is the smallest rank that a
symmetric matrix having that pattern can achieve?

Up to the addition of isolated vertices, it is easy to see that
$\G_1(F)=\{K_n\st n\in\NN\}$ for any field $F$.  In
\cite{barrett-vdHL-minrank2-infinite} and
\cite{barrett-vdHL-minrank2-finite}, $\G_2(F)$ was characterized for
any field $F$ both in terms of forbidden subgraphs and in terms of the
structure of the graph complements.  The forbidden subgraph characterizations in
these papers used ten or fewer graphs for each value of $k$.
Restricting our focus to finite fields, let $\Fld_q$ denote the finite
field with $q$ elements.  Ding and Kotlov
\cite{ding-kotlov-minrank-finite} independently used structures
similar to those introduced in this paper to obtain an upper bound for
the sizes of minimal forbidden subgraphs characterizing $\G_k(\Fld_q)$
for any $k$ and any $q$.  This result implies that there are a finite
number of forbidden subgraphs characterizing $\G_k(\Fld_q)$.  In
\cite{barrett-grout-loewy-mrF2R3}, the bound of Ding and Kotlov was
improved greatly for $\G_3(\Fld_2)$ and this set was characterized by
62 forbidden subgraphs.  This result and further computations confirm
our intuition that the forbidden subgraph characterizations of
$\G_k(\Fld_q)$ quickly become complicated as $k$ increases.

In this paper, we will characterize the structure of graphs in
$\G_k(\Fld_q)$ for any $k$ and any $q$.  The characterization is
simply stated and has a very strong connection to projective geometry
over finite fields.  At the end of the paper, we will list a few of
the ramifications of this connection to projective geometry.

We adopt the following notation dealing with fields, vector spaces,
and matrices.  Given a field $F$, the group of nonzero elements under
multiplication is denoted $F^\times$ and the vector space of dimension $k$
over $F$ is denoted $F^k$.  Given a matrix $M$, the principal
submatrix lying in the rows and columns $x_1,x_2,\ldots,x_m$ is denoted
$M[x_1,x_2,\ldots,x_m]$.

\begin{figure}[h!]
  \centering
  {\rm \begin{tikzpicture}[scale=2.54]
\ifx\dpiclw\undefined\newdimen\dpiclw\fi
\global\def\dpicdraw{\draw[line width=\dpiclw]}
\global\def\dpicstop{;}
\dpiclw=0.8bp
\dpicdraw (0.433333,0.5) circle (0.03937in)\dpicstop
\draw (0.433333,0.5) node{1};
\dpicdraw (0.1,0.166667) circle (0.03937in)\dpicstop
\draw (0.1,0.166667) node{2};
\dpicdraw (0.766667,0.166667) circle (0.03937in)\dpicstop
\draw (0.766667,0.166667) node{3};
\dpicdraw (0.1,-0.5) circle (0.03937in)\dpicstop
\draw (0.1,-0.5) node{4};
\dpicdraw (0.766667,-0.5) circle (0.03937in)\dpicstop
\draw (0.766667,-0.5) node{5};
\dpicdraw (0.362623,0.429289)
 --(0.170711,0.237377)\dpicstop
\dpicdraw (0.504044,0.429289)
 --(0.695956,0.237377)\dpicstop
\dpicdraw (0.2,0.166667)
 --(0.666667,0.166667)\dpicstop
\dpicdraw (0.1,0.066667)
 --(0.1,-0.4)\dpicstop
\dpicdraw (0.170711,0.095956)
 --(0.695956,-0.429289)\dpicstop
\dpicdraw (0.695956,0.095956)
 --(0.170711,-0.429289)\dpicstop
\dpicdraw (0.766667,0.066667)
 --(0.766667,-0.4)\dpicstop
\dpicdraw (0.2,-0.5)
 --(0.666667,-0.5)\dpicstop
\end{tikzpicture} }
  \caption{A labeled fullhouse graph}
  \label{fig:fullhouse}
\end{figure}
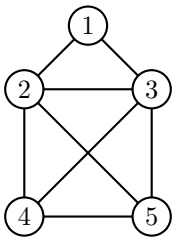

As an example of how one might approach the problem of finding the
minimum rank of a simple graph, we recall from
\cite{barrett-vdHL-minrank2-finite} the fullhouse graph in
Figure~\ref{fig:fullhouse} (there called $(P_3\cup 2K_1)^c$), which is
the only graph on 5 or fewer vertices for which the minimum rank is
field-dependent.

If $F \neq \Fld_2$, there are elements $a, b \neq 0$ in $F$ such that $a+b
\neq 0$.  Then
\begin{align*}
\begin{bmatrix}
a&a&a&0&0\\
a&a+b&a+b&b&b\\
a&a+b&a+b&b&b\\
0&b&b&b&b\\
0&b&b&b&b
\end{bmatrix} \in S(F,\fh)
\end{align*}
which shows that $\mr(F,\fh)=2$.  The case $F=\Fld_2$ gives a different
result.  Let $A$ be any matrix in $S(\Fld_2,\fh)$.  Then for some
$d_1,d_2,\ldots, d_5\in \Fld_2$,
\begin{align*}
A=  \begin{bmatrix}
d_1&1&1&0&0\\
1&d_2&1&1&1\\
1&1&d_3&1&1\\
0&1&1&d_4&1\\
0&1&1&1&d_5
\end{bmatrix}
\quad \text{ and } \quad \det (A[\{1, 2,
5\},\{1, 3, 4\}])=\left| \begin{array}{ccc}
    d_1&1&0\\1&1&1\\0&1&1\end{array}\right| =1,
\end{align*}
where $A[\{1, 2, 5\},\{1, 3, 4\}]$ is the submatrix of $A$ lying in rows
$\{1,2,5\}$ and columns $\{1,3,4\}$.  Therefore $\mr(\Fld_2,\fh)\geq 3$.
Setting each $d_i$ to 1 verifies that $\mr(\Fld_2,\fh)= 3$.

In spite of this dependence on the field, there are a number of
results about minimum rank that are field independent.  For example,
the minimum rank of a tree is field independent (see any of
\cite{bank-thesis}, \cite{sinkovic-thesis}, or
\cite{hogben-minrank-tree}).  Many of the forbidden subgraphs
classifying $\G_3(\Fld_2)$ that are found in
\cite{barrett-grout-loewy-mrF2R3} are also forbidden subgraphs for
$\G_3(F)$ for any field $F$.  These results and others demonstrate
that results obtained over finite fields can provide important
insights for other fields.

The presentation of material in this paper is oriented towards a
reader that is familiar with concepts from linear algebra and graph
theory.  In the rest of this section, we will review some of our
conventions in terminology from graph theory.

In this paper, graphs are undirected, may have loops, but will not
have multiple edges between vertices.  To simplify our drawings, a
vertex with a loop (a \emph{looped vertex}) will be filled (black) and
a vertex without a loop (a \emph{nonlooped vertex}) will be empty
(white).  A \emph{simple graph} is a graph without loops.  Let $G$ be
a graph with some loops and $\hat G$ be the simple version of $G$ obtained
by deleting all loops.  We say that a matrix in $S(F,\hat G)$
\emph{corresponds} to the simple graph $\hat G$.  A matrix $A\in S(F,\hat G)$
\emph{corresponds} to $G$ if $a_{ii}$ is nonzero exactly when the
vertex $i$ has a loop in $G$.  Note that if a matrix corresponds to a
looped graph, then it also corresponds to the simple version of the
graph.

We recall some notation from graph theory.

\begin{definition} 
  Given two graphs $G$ and $H$ with disjoint vertex sets $V(G)$ and
  $V(H)$ and edge sets $E(G)$ and $E(H)$, the \emph{union} of $G$ and
  $H$, denoted $G \cup H$, has vertices $V(G)\cup V(H)$ and edges $E(G)\cup
  E(H)$.  The \emph{join} of $G$ and $H$, denoted $G\lor H$, has vertices
  $V(G)\cup V(H)$ and edges $E(G)\cup E(H)\cup\{uv \st u\in V(G),\, v\in V(H)\}$.
  The complement of the graph $G$, denoted $G^c$, has vertices $V(G)$
  and edges $\{uv \st u,v\in V(G),\, uv\not\in E(G)\}$.  Note that a vertex
  is looped in $G$ if and only if it is nonlooped in $G^c$.
\end{definition}

\begin{definition} 
  The simple complete graph on $n$ vertices will be denoted by $K_n$
  and has vertices $\{1,2,\ldots,n\}$ and edges $\{xy\st x,y\in V(K_n), x\neq y\}$.
  The simple complete multipartite graph $K_{s_1,s_2,\ldots,s_m}$ is
  defined as $K_{s_1}^c\lor K_{s_2}^c\lor\cdots \lor K_{s_m}^c$.
\end{definition}

\begin{definition} 
  Two vertices in a graph are \emph{adjacent} if an edge connects them.  A
  {\it clique} in a graph is a set of pairwise adjacent vertices.  An
  {\it independent set} in a graph is a set of pairwise nonadjacent
  vertices.
\end{definition}

The next definition extends a standard definition introduced in
\cite{blowup-lemma} and is used in random graph theory in connection
with the regularity lemma.

\begin{definition}
  A \emph{blowup} of a graph $G$ with vertices $\{v_1,v_2,\ldots,v_n\}$ is a
  new simple graph $H$ constructed by replacing each nonlooped vertex
  $v_i$ in $G$ with a (possibly empty) independent set $V_i$, each
  looped vertex $v_i$ with a (possibly empty) clique $V_i$, and each
  edge $v_iv_j$ in $G$ ($i\neq j$) with the edges $\{xy \st x\in V_i, y\in
  V_j\}$ in $H$.
\end{definition}

\begin{example}\label{ex:blowup}
Let $G$ be the graph labeled in Figure~\ref{fig:example-blowup}\subref{fig:example-blowup-orig}.
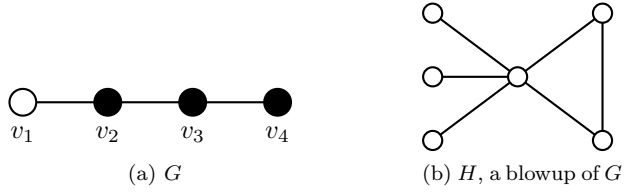
\begin{figure}[h!]
  \centering
 \subfloat[$G$]{\begin{tikzpicture}[scale=2.54]
\ifx\dpiclw\undefined\newdimen\dpiclw\fi
\global\def\dpicdraw{\draw[line width=\dpiclw]}
\global\def\dpicstop{;}
\dpiclw=0.8bp
\dpicdraw (0.07,0.05) circle (0.027559in)\dpicstop
\dpicdraw[fill=black](0.514444,0.05) circle (0.027559in)\dpicstop
\dpicdraw[fill=black](0.958889,0.05) circle (0.027559in)\dpicstop
\dpicdraw[fill=black](1.403333,0.05) circle (0.027559in)\dpicstop
\dpicdraw (0.14,0.05)
 --(0.444444,0.05)\dpicstop
\dpicdraw (0.584444,0.05)
 --(0.888889,0.05)\dpicstop
\dpicdraw (1.028889,0.05)
 --(1.333333,0.05)\dpicstop
\draw (0.07,-0.12) node{$v_1$};
\draw (0.514444,-0.12) node{$v_2$};
\draw (0.958889,-0.12) node{$v_3$};
\draw (1.403333,-0.12) node{$v_4$};
\end{tikzpicture}\label{fig:example-blowup-orig}}
 \qquad\qquad
 \subfloat[$H$, a blowup of~$G$]{\begin{tikzpicture}[scale=2.54]
\ifx\dpiclw\undefined\newdimen\dpiclw\fi
\global\def\dpicdraw{\draw[line width=\dpiclw]}
\global\def\dpicstop{;}
\dpiclw=0.8bp
\dpicdraw (0.05,-0.333333) circle (0.019685in)\dpicstop
\dpicdraw (0.05,0) circle (0.019685in)\dpicstop
\dpicdraw (0.05,0.333333) circle (0.019685in)\dpicstop
\dpicdraw (0.494444,0) circle (0.019685in)\dpicstop
\dpicdraw (0.938889,-0.333333) circle (0.019685in)\dpicstop
\dpicdraw (0.938889,0.333333) circle (0.019685in)\dpicstop
\dpicdraw (0.09,-0.303333)
 --(0.454444,-0.03)\dpicstop
\dpicdraw (0.1,0)
 --(0.444444,0)\dpicstop
\dpicdraw (0.09,0.303333)
 --(0.454444,0.03)\dpicstop
\dpicdraw (0.534444,-0.03)
 --(0.898889,-0.303333)\dpicstop
\dpicdraw (0.534444,0.03)
 --(0.898889,0.303333)\dpicstop
\dpicdraw (0.938889,-0.283333)
 --(0.938889,0.283333)\dpicstop
\end{tikzpicture}\label{fig:example-blowup-blowup}}
 \caption{Graphs in Example~\ref{ex:blowup}}
   \label{fig:example-blowup}
\end{figure}

Let $|V_1|=3$, $|V_2|=1$, $|V_3|=2$, and $|V_4|=0$.  Then we obtain
the simple blowup graph $H$ in
Figure~\ref{fig:example-blowup}\subref{fig:example-blowup-blowup}.  It
is useful to see how matrices corresponding to a graph and a blowup of
the graph are related.  Over $\Fld_3$, let
\begin{align*}
  M=\left[ \begin{array}{c|c|c|c}
        0&2&0&0\\ \hline
        2&1&1&0\\ \hline
        0&1&1&1\\ \hline
        0&0&1&1 \end{array}\right]
  \quad \text{ and } \quad
  N=\left[ \begin{array}{ccc|c|cc}
 0&0&0&1&0&0\\ 0&0&0&2&0&0\\ 0&0&0&1&0&0\\\hline
 1&2&1&0&1&1\\\hline
  0&0&0&1&0&1\\  0&0&0&1&1&2 \end{array}\right].
\end{align*}
Then $M$ is an example of a matrix corresponding to $G$ and $N$ is an
example of a matrix corresponding to $H$. Note that, for example, the
entry $m_{11}$ was replaced with a $3\times 3$ zero block in $N$, the entry
$m_{12}$ was replaced with a $3\times 1$ nonzero block in $N$, the entries
in the last row and column of $M$ were replaced with empty blocks
(i.e., erased), and the diagonal entries of $N$ were changed to
whatever was desired.  These substitutions of block matrices
correspond to the vertex substitutions used to construct $H$.
\end{example}

We will introduce our method by presenting a proof of a special case
of a characterization theorem from \cite{barrett-vdHL-minrank2-finite}
which characterizes $\G_2(\Fld_2)$.  We will then generalize this
proof into a characterization of all simple graphs in $\G_{k}(\Fld_q)$
for any $k$ and $q$.  After giving examples for some specific $k$ and
$q$, we will describe the strong connection to projective geometry and
list some consequences of this connection.

\section{A new approach to a recent result}

We will introduce our method by giving a proof of a special case of
Theorems~5 and 6 of \cite{barrett-vdHL-minrank2-finite}.

\begin{theorem}[{\cite{barrett-vdHL-minrank2-finite}}]
  \label{thm:f2r2-orig}
  Let $G$ be a simple graph on $n$ vertices.  Then $\mr(\Fld_2,G)\leq2$ if
  and only if the simple version of $G^c$ is either of the form
  \begin{equation*}
    (K_{s_1}\cup K_{p_1,q_1})\lor K_r
  \end{equation*}
  for some appropriate nonnegative integers $s_1$, $p_1$, $q_1$, and
  $r$, or of the form
  \begin{equation*}
    (K_{s_1}\cup K_{s_2}\cup K_{s_3})\lor K_r
  \end{equation*}
  for some appropriate nonnegative integers $s_1$, $s_2$, $s_3$, and $r$.
\end{theorem}

We first rephrase Theorem~\ref{thm:f2r2-orig} using blowup graph terminology.

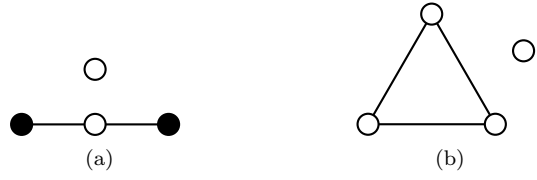
\begin{figure}[h!]
  \centering
 \subfloat[]{\begin{tikzpicture}[scale=2.54]
\ifx\dpiclw\undefined\newdimen\dpiclw\fi
\global\def\dpicdraw{\draw[line width=\dpiclw]}
\global\def\dpicstop{;}
\dpiclw=0.8bp
\dpicdraw[fill=black](0.055,-0.144231) circle (0.021654in)\dpicstop
\dpicdraw (0.439615,-0.144231) circle (0.021654in)\dpicstop
\dpicdraw[fill=black](0.824231,-0.144231) circle (0.021654in)\dpicstop
\dpicdraw (0.439615,0.144231) circle (0.021654in)\dpicstop
\dpicdraw (0.11,-0.144231)
 --(0.384615,-0.144231)\dpicstop
\dpicdraw (0.494615,-0.144231)
 --(0.769231,-0.144231)\dpicstop
\end{tikzpicture}\label{fig:thm-f2r2-a}}
 \qquad\qquad\qquad
 \subfloat[]{\begin{tikzpicture}[scale=2.54]
\ifx\dpiclw\undefined\newdimen\dpiclw\fi
\global\def\dpicdraw{\draw[line width=\dpiclw]}
\global\def\dpicstop{;}
\dpiclw=0.8bp
\dpicdraw (0.388087,0.288462) circle (0.021654in)\dpicstop
\dpicdraw (0.055,-0.288461) circle (0.021654in)\dpicstop
\dpicdraw (0.721173,-0.288462) circle (0.021654in)\dpicstop
\dpicdraw (0.868856,0.096154) circle (0.021654in)\dpicstop
\dpicdraw (0.360587,0.24083)
 --(0.0825,-0.24083)\dpicstop
\dpicdraw (0.11,-0.288461)
 --(0.666173,-0.288462)\dpicstop
\dpicdraw (0.693673,-0.24083)
 --(0.415587,0.24083)\dpicstop
\end{tikzpicture}\label{fig:thm-f2r2-b}}

 \caption{Graphs in Theorem~\ref{thm:f2r2}}
   \label{fig:thm-f2r2}

\end{figure}

\begin{theorem}[{\cite{barrett-vdHL-minrank2-finite}}] 
  \label{thm:f2r2}
  Let $G$ be a simple graph on $n$ vertices.  Then
  \mbox{$\mr(\Fld_2,G)\leq2$} (i.e., $G\in \G_2(\Fld_2)$) if and only if
  $G$ is a blowup of either of the graphs in Figure~\ref{fig:thm-f2r2}.
\end{theorem}

In the proof of this result, we will need the following lemma and
corollary, which hold in any field.  We will then give a proof of Theorem~\ref{thm:f2r2}.

\begin{lemma}[{\cite[Theorem~8.9.1]{gr-graph-theory}}] 
\label{lem:gr-rank-decomp}
Let $A$ be an $n\times n$ symmetric matrix of rank~$k$.  Then there
is an invertible principal $k\times k$ submatrix $B$ of $A$ and a
$k\times n$ matrix $U$ such that
\begin{align*}
  A=U^t B U.
\end{align*}
\end{lemma}

\begin{corollary}
\label{cor:rank-decomp}
Let $A$ be an $n\times n$ symmetric matrix.  Then $\rank A \leq k$ if
and only if there is some invertible $k\times k$ matrix $B$ and
$k\times n$ matrix $U$ such that $A=U^tBU$.
\end{corollary}

\begin{proof}
  Let $A$ have rank $r\leq k$.  Then by
  Lemma~\ref{lem:gr-rank-decomp}, there is an invertible $r\times r$
  matrix $B_1$ and an $r\times n$ matrix $U_1$ such that
  $A=U_1^tB_1U_1$.  Let $B_2=\bmat B_1&O\\O&I_{k-r}\emat$ and
  $U_2=\bmat U_1\\O \emat$ (where $O$ represents a zero matrix of the
  appropriate size).  Then $A=U_2^tB_2U_2$.

  The reverse implication follows from the rank inequality
  $\rank(U^tBU)\leq \rank B$.
\end{proof}

Recall that two square matrices $A$ and $B$ are congruent if there
exists some invertible matrix $C$ such that $A=C^tBC$.  It is
straightforward to show that congruence is an equivalence relation.
Let $\mathcal{B}$ consist of one representative from each congruence
equivalence class of invertible symmetric $k\times k$ matrices.  By
Corollary~\ref{cor:rank-decomp}, if $A$ is a symmetric $n\times n$ matrix
with $\rank A\leq k$, then $A\in \{U^tBU \st B\in\mathcal{B},\, U \text{ a $k\times
  n$ matrix}\}$.

We now proceed with the proof of Theorem~\ref{thm:f2r2}.  

\begin{proof}[Proof of Theorem~\ref{thm:f2r2}]
  First, we compute a suitable $\mathcal{B}$, a set of representatives
  from the congruence classes of invertible symmetric $2\times 2$ matrices
  over $\Fld_2$.  If an invertible symmetric $2\times 2$ matrix $B$ over
  $\Fld_2$ has a nonzero diagonal entry, then $B=\bmat 1&1\\1&0
  \emat$, $B=\bmat 0&1\\1&1\emat$, or $B=I_2$.  In any of these three
  cases, $B^tBB=I_2$, so $B$ is congruent to the identity matrix
  $I_2$.  If an invertible symmetric $2\times 2$ matrix $B$ over $\Fld_2$
  has all zeros on the diagonal, then the off-diagonal entries must be
  nonzero, so $B=\bmat 0&1\\1&0\emat$.  In this case,
\begin{align*}
\bmat a&c\\b&d\emat \bmat 0&1\\1&0\emat \bmat a&b\\c&d\emat = \bmat ac+ac&ad+bc\\ad+bc&bd+bd\emat=
\bmat 0&ad+bc\\ad+bc&0\emat,
\end{align*}
so any matrix congruent to $B$ will have a zero diagonal.
Therefore, a suitable $\mathcal{B}$ is
\begin{align*}
  \mathcal{B}=\left\{I_2,\bmat 0&1\\1&0\emat\right\}.
\end{align*}

Because $U$ is a matrix with entries in $\Fld_2$, the columns of $U$
are members of the finite set
\begin{align*}
\left\{ \bmat1\\0 \emat, \bmat 0\\1 \emat, \bmat 1\\1 \emat, 
\bmat 0\\0\emat \right\}.
\end{align*}
Let $A$ be a symmetric $k\times k$ matrix.  For any $n \times n$
permutation matrix $P$, the graphs of $A$ and $P^tAP$ are isomorphic.
Therefore we may assume that identical columns of $U$ are contiguous
and write $U = \bmat E_1 & E_2 & J & O \emat$ where $E_1$ is $2 \times
p$ matrix with each column equal to $\bmat 1\\0\emat$, $E_2$ is $2
\times q$ matrix with each column equal to $\bmat 0\\1\emat$, $J$ is a
$2 \times r$ matrix with each entry equal to $1$, and $O$ is a $2
\times t$ zero matrix. Then either
$$A=\bmat
E_1^{\rm{T}}\\
E_2^{\rm{T}}\\
J^{\rm{T}}\\
O^{\rm{T}}\emat \bmat
E_1 & E_2 & J & O \emat = \bmat
J_p & O & J_{p,r} & O\\
O & J_q & J_{q,r} & O\\
J_{r,p} & J_{r,q} & O_r & O\\
O & O & O & O_t\emat$$
or else $$A=\bmat
E_1^{\rm{T}}\\
E_2^{\rm{T}}\\
J^{\rm{T}}\\
O^{\rm{T}}\emat \bmat
0 & 1\\
1 & 0\emat \bmat
E_1 & E_2 & J & O \emat = \bmat
O_p &  J_{p,q} & J_{p,r} & O\\
J_{q,p} & O_q & J_{q,r} & O\\
J_{r,p} & J_{r,q} & O_r & O\\
O & O & O & O_t\emat,$$
where $J$ is an all-ones matrix, $O$ is a zero matrix, and subscripts of $J$ and $O$ denote the dimensions of the matrix.

Any simple graph corresponding to the first matrix is a blowup of the
graph in Figure~\ref{fig:thm-f2r2}\subref{fig:thm-f2r2-a}, while any
simple graph corresponding to the second matrix is a blowup of the
graph in Figure~\ref{fig:thm-f2r2}\subref{fig:thm-f2r2-b}.  Thus we
have established Theorem~\ref{thm:f2r2}.
\end{proof}

\begin{observation} \label{obs:simplified-gFk}  
Note that every block in the above matrices is either a $O$ matrix or
a $J$ matrix.  Consequently, we could have obtained the zero/nonzero
form of the matrices with rank at most 2 by only considering $U = \bmat
1 & 0 & 1 & 0\\
0 & 1 & 1 & 0\emat$ and computing
\begin{align*}
A = U^tU = \bmat
1 & 0 & 1 & 0\\
0 & 1 & 1 & 0\\
1 & 1 & 0 & 0\\
0 & 0 & 0 & 0\emat
\end{align*}
and
\begin{align*}
A = U^tB_2U = \bmat
1 & 0\\
0 & 1\\
1 & 1\\
0 & 0\emat \bmat
0 & 1 & 1 & 0\\
1 & 0 & 1 & 0\emat = \bmat
0 & 1 & 1 & 0\\
1 & 0 & 1 & 0\\
1 & 1 & 0 & 0\\
0 & 0 & 0 & 0\emat.
\end{align*}
The nonzero diagonal entries correspond to loops in our graphs.  This
simplified procedure again yields the graphs in
Figure~\ref{fig:thm-f2r2}.
\end{observation}

In the proof of Theorem~\ref{thm:f2r2}, we noted that any $U$ could be
written in a standard form.  In Observation~\ref{obs:simplified-gFk},
we saw how the standard form of $U$ could be simplified to take advantage
of the theorem being about blowup graphs.  We will now discuss the
reasoning behind these constructions and show that an analogous
standard form of $U$ exists for any finite field and any $k$.

Because we construct the graphs using representatives of
congruence classes, it is important for any simplified $U$ to have the
property that if $B$ and $\hat B$ are congruent, then $U^tBU$ and
$U^t\hat B U$ correspond to isomorphic graphs.  The following lemma
shows that if we take a matrix $U$ where the columns consist of all
vectors in $\Fld_q^k$, like in Observation~\ref{obs:simplified-gFk},
and if $B$ and $\hat B$ are congruent, then $U^tBU$ and $U^t\hat BU$
correspond to isomorphic graphs.

\begin{lemma}\label{lem:U-all-vectors}
  Let $U$ be the matrix with columns $\{v \st v\in \Fld_q^k\}$.
  Let $B$ and $C$ be invertible $k\times k$ matrices with $B$ symmetric.  Then the graphs corresponding to
  $U^tBU$ and $U^t(C^tBC)U$ are isomorphic.
\end{lemma}

\begin{proof}
  Since every vector in $\Fld_q^k$ appears as a column of $U$ and the mapping
  $x\mapsto Cx$ is one-to-one, $CU$ is just a column permutation of
  $U$.  This permutation induces a relabeling of the graph $U^tBU$ to
  give the graph of $(CU)^tB(CU)=U^t(C^tBC)U$.
\end{proof}

Though this invariance property with respect to congruent matrices does not
hold for an arbitrary $U$, there is another smaller $U$ which does
have the same property.  We first need some preliminary material.
Then we will introduce this new $U$ in Lemma~~\ref{lem:congruent-isomorphic}.

\begin{definition}
  Let $F$ be a field.  Two nonzero vectors $v_1,v_2\in F^k$ are
  \emph{projectively equivalent} if there exists some nonzero $c\in F$
  such that $v_1=cv_2$.
\end{definition}

It is easy to check that projective equivalence is in fact
an equivalence relation on the vectors in $V$.

We pause to note that replacing a column of $U$ with a projectively
equivalent column does not affect the graph corresponding to $U^tBU$.
To see this, let $U=[u_1\; u_2\; \cdots \; u_n]$ and let $i\in\{1,2,\ldots,n\}$.
Let $\hat U$ be the matrix obtained from $U$ by replacing the column
$u_i$ with $cu_i$ for some nonzero $c\in F$.  Then the $i,j$ entry of
$\hat U^tB\hat U$, $(cu_i)^tBu_j$ if $i\neq j$ or $(cu_i)^tB(cu_i)$ if
$i=j$, is zero if and only if the $i,j$ entry of $U^tBU$, $u_i^tBu_j$,
is zero.  Thus the graphs associated with $U^tBU$ and $\hat U^t B \hat
U$ are equal.

\begin{lemma}\label{lem:projective-permutation}
  Let $F$ be any field, let $x\in F^k$, let $\x$ denote
  the projective equivalence class of $x$, and let $P=\cup_{x\in
    F^k-\0}\{\x\}$, the set of projective equivalence classes in $F^k$.
  Let $C$ be an invertible matrix.  Then the map $f\colon P \to P$
  defined by $f\colon \x \mapsto \class{Cx}$ is a bijection.
\end{lemma}

\begin{proof}
  The function $f$ is well-defined since if $Cx=y$, then for any
  nonzero $k\in F$, $\class{C(kx)}=\class{kCx}=\class{ky}=\y$.  If $
  \class{Cx_1}=\class{Cx_2}$, then for some nonzero $k\in F$,
  $kCx_1=Cx_2$, which implies $C(kx_1-x_2)=0$, giving $kx_1=x_2$ since
  $C$ is invertible.  Therefore $\class{x_1}=\class{x_2}$ and $f$ is
  injective.  Surjectivity of $f$ also follows from the hypothesis
  that $C$ is invertible.
\end{proof}

\begin{lemma}\label{lem:congruent-isomorphic}
  Let $\class{x_1}, \class{x_2},\ldots,\class{x_m}$ be the projective
  equivalence classes of $\Fld_q^k-\0$, with each $x_i$ as a chosen
  representative from its class.  Let $U=[x_1\;x_2\;\cdots\; x_m]$, the
  matrix with column vectors $x_1,x_2,\ldots,x_m$.  Let $B$ and $C$ be
  invertible $k\times k$ matrices with $B$ symmetric.  Then the graphs
  corresponding to $U^tBU$ and $U^t(C^tBC)U$ are isomorphic.
\end{lemma}

\begin{proof}
  Let $T=CU$.  Denote the $i$th column of $U$ by $u_i$ and the $i$th
  column of $T$ by $t_i$.  By Lemma~\ref{lem:projective-permutation},
  the sequence of projective equivalence classes $\class{t_1},
  \class{t_2}, \ldots, \class{t_n}$ is just a permutation of the sequence
  $\class{u_1},\class{u_2},\ldots,\class{u_n}$.  Form the matrix $S$ in
  which the $i$th column, $s_i$, is $u_j$ if
  $\class{t_i}=\class{u_j}$, so that $S$ is a column permutation of
  $U$ and $\class{s_i}=\class{t_i}$.  Then the graph corresponding to
  $U^t(C^tBC)U=(C U)^tB(C U)=T^tBT$ is isomorphic to the graph
  corresponding to $S^tBS$ by the reasoning preceding
  Lemma~\ref{lem:projective-permutation}, which is in turn just a
  relabeling of the graph corresponding to $U^tBU$.
\end{proof}

We now find a standard form for any matrix $U$, as in our proof of
Theorem~\ref{thm:f2r2}.  Let $U$ be a $k\times n$ matrix over $\Fld_q$ and
let $B$ be an invertible symmetric $k\times k$ matrix over $\Fld_q$.  Let
$\class{x_1}, \class{x_2},\ldots,\class{x_m}$ be the projective equivalence
classes of $\Fld_q^k-\0$, with each $x_i$ as a chosen representative
from its class.  For each nonzero column $u_i$ of $U$, replace $u_i$
with the chosen representative of $\class{u_i}$.  Then permute the
columns of $U$ so that the matrix is of the form $\hat U=[X_1 \; X_2
\; \cdots\; X_m\; O]$, where each $X_i$ is a block matrix of columns equal
to $x_i$ and $O$ is a zero block matrix.  Note that some of these
blocks may be empty.  Let $G$ be the simple graph corresponding to
$U^tBU$ and let $\hat G$ be the simple graph corresponding to $\hat
U^tB\hat U$.  From our results above, $G$ is isomorphic to $\hat G$.

As illustrated in Observation~\ref{obs:simplified-gFk}, we can obtain
the zero/nonzero structure of the block matrix $\hat U^tB\hat U$ by
simply deleting all duplicate columns of $\hat U$.  Deleting these
duplicate columns of $\hat U$ leaves a matrix that can be obtained
from $\tilde U=[x_1\;x_2\;\cdots\;x_m\;0]$ by deleting the columns of
$\tilde U$ corresponding to empty blocks of $\hat U$.  Let $\tilde G$
be the (looped) graph corresponding to $\tilde U^tB\tilde U$.  Then
$\hat G$ is a blowup of $\tilde G$, which implies that $G$ is a blowup
of $\tilde G$.

Furthermore, let $\mathcal{B}$ be a set consisting of one
representative from each congruence class of invertible symmetric $k\times
k$ matrices and let $\hat B$ be the representative that is congruent
to $B$.  Then from Lemma~\ref{lem:congruent-isomorphic}, the graphs
corresponding to $\tilde U^t B \tilde U$ and $\tilde U^t\hat B \tilde
U$ are isomorphic.

There is another simplification we can make.  Notice that both graphs
displayed in Theorem~\ref{thm:f2r2} have an isolated nonlooped vertex.
This vertex came from the zero column vectors in $U$ and corresponds
to the fact that adding any number of isolated vertices to a graph
does not change its minimum rank.  In any theorem like
Theorem~\ref{thm:f2r2}, each graph from which we construct blowups
will always have this isolated nonlooped vertex and so will be of the
form $G\cup K_1$.  Note that in constructing such a graph $G$, it is
enough to assume that $\tilde U$ in the above paragraphs does not have
a zero column vector.

\begin{definition}\label{def:g_k}
  Let $\class{x_1}, \class{x_2},\ldots,\class{x_m}$ be the projective
  equivalence classes of $\Fld_q^k-\0$, with each $x_i$ as a chosen
  representative from its class.  Let $\mathcal{B}$ be a set
  consisting of one representative from each congruence class of
  invertible symmetric $k\times k$ matrices.  Let $U=[x_1\;x_2\;\cdots\; x_m]$,
  the matrix with column vectors $x_1,x_2,\ldots,x_m$.  We define the set
  of graphs $\g_k(\Fld_q)$ as the set of graphs corresponding to the
  matrices in $\{U^tBU \st B\in \mathcal{B}\}$.
\end{definition}

We now have the following result (recall that $K_1$ has no loop).

\begin{theorem}
  A simple graph $G$ is in $\G_k(\Fld_q)$ if and only if $G$ is a
  blowup of some graph in $\{H\cup K_1 \st H\in\g_k(\Fld_q)\}$.
\end{theorem}

\begin{proof}
  Let $G$ be a simple graph in $\G_k(\Fld_q)$.  Let $A\in S(\Fld_q,G)$ be
  a matrix with $\rank A\leq k$.  Then $A=U^tBU$ for some $k\times n$ matrix
  $U$ and some invertible symmetric $k\times k$ matrix $B$.  Using the
  procedure outlined in the paragraphs following
  Lemma~\ref{lem:congruent-isomorphic}, we see that $G$ is a blowup of
  a graph $\tilde G$ corresponding to $\tilde U^t B \tilde U$, where
  $\tilde U$ and $B$ are defined as in the procedure.
  Lemma~\ref{lem:congruent-isomorphic} then shows that $\tilde
  G\in\g_k(\Fld_q)$.

  Conversely, let $G$ be a blowup of some graph in $\{H\cup K_1\st
  H\in\g_k(\Fld_q)\}$ obtained by replacing each vertex $v_i$ of $H$ with
  a set of vertices $V_i$ and $K_1$ with any number of vertices.
  Deleting isolated vertices of $G$ does not change the minimum rank
  of $G$, so without loss of generality, we will assume that $G$ has
  no isolated vertices (which implies that $K_1$ was replaced with an
  empty set of vertices).  Let $\class{x_1},
  \class{x_2},\ldots,\class{x_m}$ be the projective equivalence
  classes of $\Fld_q^k-\0$, with each $x_i$ as a chosen representative
  from its class.  Let $\tilde U= [x_1\; x_2\; \cdots \; x_m]$ and let $B$
  be an invertible symmetric $k\times k$ matrix such that $\tilde U^t B
  \tilde U$ corresponds to the graph $H$.  Form the matrix $\hat
  U=[X_1\;X_2\;\cdots \; X_m]$ by replacing each column $x_i$ of $\tilde U$
  with the block $X_i$, where the columns of $X_i$ consist of $|V_i|$
  copies of $x_i$.  Then $\hat U^t B \hat U$ corresponds to $G$ and
  $\rank \hat U^t B \hat U\leq k$ since $B$ has rank $k$.  Thus
  $\mr(\Fld_q,G)\leq k$, so $G\in\G_k(\Fld_q)$.
\end{proof}

Now we will make this into a more explicit characterization of
$\G_k(\Fld_q)$ by finding a suitable $\mathcal{B}$ for any $k$ and any
$q$, thus enabling us to explicitly find $\g_k(\Fld_q)$ for any $k$
and any $q$.

\section{Congruence classes of symmetric matrices over finite fields}
\label{sec:congr-class-symm}
Symmetric matrices represent symmetric bilinear forms and play an
important role in projective geometry.  Two congruent symmetric
matrices represent the same symmetric bilinear form with respect to
different bases.  Because of their fundamental importance, congruence
classes of symmetric matrices over finite fields have been studied and
characterized for a long time in projective geometry.  In this
section, we have distilled the pertinent proofs of these
characterizations from \cite{albert-matrices},
\cite{hirschfeld-projective}, and \cite{cohn-algebra} to give a
suitable $\mathcal{B}$ for invertible symmetric $k\times k$ matrices
over $\Fld_q$ for any $k$ and $q$.  In the next section, we will
expound more on the connection between the minimum rank problem and
projective geometry.

We need the following elementary lemma.

\begin{lemma}\label{lem:make-diagonal}
  If a symmetric matrix $B=\bmat C &D\\D^t&E \emat$, where $C$ is a
  square invertible matrix, then $B$ is congruent to $\bmat C&O\\O&E'
  \emat$, where $O$ is a zero matrix and $E'$ is a square symmetric
  matrix of the same order as $E$.
\end{lemma}

\begin{proof}
  Let $R=C^{-1}D$ so that $CR=D$.  Then
  \begin{align*}
    \bmat I&O\\-R^t&I \emat \bmat C&D\\D^t&E\emat \bmat I&-R\\O&I\emat
    &= \bmat C&D\\-R^tC+D^t&-R^tD+E\emat \bmat I&-R\\O&I\emat \\
    &= \bmat C&-CR+D\\-R^tC+D^t&R^tCR-D^tR-R^tD+E\emat\\
    &= \bmat C&O\\O&E-D^tR\emat,
  \end{align*}
  since $-CR+D=O=(-CR+D)^t=-R^tC+D^t$.
\end{proof}

\begin{lemma}\label{lem:canonical-matrix}
  Every symmetric matrix over $\Fld_q$ is congruent to a matrix of
  the form $\diag(a_1,a_2,\ldots, a_s,b_1H_1,b_2H_2,\ldots, b_tH_t)$, where
  $a_i,b_i\in \Fld_q$, $H_i=\bmat 0&1\\1&0\emat$, and $s$ and $t$ are
  nonnegative integers.
\end{lemma}

\begin{proof}
If $B$ is the zero matrix, then the result is true.

If $B$ is not the zero matrix, then the diagonal of $B$ has a
nonzero entry or there is some $a_{ij}\neq 0$, $i\neq j$, so that $B$
has a principal submatrix of the form $\bmat 0&a_{ij}\\a_{ij} &0
\emat=a_{ij}H$, where $H=\bmat 0&1\\1&0\emat$.  

In the first case, by using a suitable permutation, we may assume that
$b_{11}\neq 0$.  By Lemma~\ref{lem:make-diagonal}, $B$ is congruent to
$\diag(b_{11},B')$.

In the second case, again by using a suitable permutation, we may assume that 
the upper left $2\times 2$ principal submatrix is $a_{ij}H$.  By
Lemma~\ref{lem:make-diagonal}, $B$ is congruent to $\diag(a_{ij}H,B')$.

Continue this process inductively with $B'$.  Then, again using a
suitable permutation, $B$ is congruent to $\diag(a_1,a_2,\ldots,
a_s,b_1H,b_2H,\ldots, b_tH)$.
\end{proof}

We will now treat the even characteristic and odd characteristic cases
separately.

\subsection{Even characteristic}
We first consider the case when $\Fld_q$ has even characteristic.
First, we need a well-known result.
\begin{lemma}
  Every element in a field of characteristic 2 is a square.
\end{lemma}

\begin{corollary}\label{cor:diag-odd-char}
  Every symmetric matrix is congruent to $\diag(I_s,H_1,H_2,\ldots, H_t)$.
\end{corollary}

\begin{proof}
  By Lemma~\ref{lem:canonical-matrix}, a symmetric matrix $A$ is
  congruent to a matrix 
  \begin{align*} 
    B=\diag(a_1,a_2,\ldots, a_s,b_1H_1,b_2H_2,\ldots, b_tH_t).
  \end{align*}
  Let
  \begin{align*}  
    C=\diag(\frac{1}{\sqrt{a_1}},\frac 1 {\sqrt{a_2}},\ldots,\frac 1
    {\sqrt{a_s}},\frac 1 {\sqrt{b_1}}I_2,\frac 1
    {\sqrt{b_2}}I_2,\ldots,\frac 1 {\sqrt{b_t}}I_2).
  \end{align*}
  Then $C^tBC=\diag(I_s,H_1,H_2,\ldots,H_t)$.
\end{proof}

Let $B$ be a symmetric matrix in $\Fld_q$.  Then according to
Corollary~\ref{cor:diag-odd-char}, $B$ is congruent to a matrix
$C=\diag(I_s,H_1,H_2,\ldots,H_t)$, where each
$H_i=[\begin{smallmatrix}0&1\\1&0\end{smallmatrix}]$.  Either $s=0$ or $s>0$.
If $s>0$, then $\diag(I_s,H_1,H_2,\ldots,H_t)$, and thus $B$, is congruent
to $I_k$.  To see this, let
\begin{align*}
  A=\diag(1,H)=\bmat 1&0&0\\0&0&1\\0&1&0 \emat \quad \text{ and } \quad
  C=\bmat 1&1&1\\1&0&1\\0&1&1 \emat.
\end{align*}
Then, since $\characteristic \Fld_q=2$,
\begin{align*}
  C^t(AC)=\bmat 1&1&0\\1&0&1\\1&1&1 \emat \bmat 1&1&1\\0&1&1\\1&0&1
  \emat  = I_3.
\end{align*}
If $s=0$, then $\diag(H_1,H_2,\ldots,H_t)$ and $B$ have even order and $B$
is congruent to $\diag(H_1,\ldots,H_{k/2})$.

The next lemma shows that these two cases are different.

\begin{lemma}
  If a symmetric matrix $B$ has a zero diagonal, then every matrix
  congruent to $B$ has a zero diagonal.
\end{lemma}

\begin{proof}
  Let $B$ be a symmetric matrix having a zero diagonal.  If $v$ is the
  $k$th column of a matrix $C$, then the $(k,k)$ entry of $C^tBC$ is
  $v^tBv$, which is zero, since
\begin{equation*}
  v^tBv=\sum_{i,j}b_{ij}v_iv_j=\sum_{i}b_{ii}v_i^2+\sum_{i<j}b_{ij}(v_iv_j+v_iv_j)
  = \sum_ib_{ii}v_i^2=0.\qedhere
\end{equation*}
\end{proof}

The results in this subsection give us the following lemma.

\begin{lemma}
  \label{lem:B-for-even-char}
  Let $q$ be even.  To determine $\g_k(\Fld_q)$, we may take $\mathcal{B}$
  as follows: if $k$ is odd, then $\mathcal{B}=\{I_k\}$; if $k$ is
  even, then $\mathcal{B}=\{I_k,\diag(H_1,H_2,\dots,H_{k/2})\}$, where
  $H_i=\bmat 0&1\\1&0 \emat$.
\end{lemma}

\subsection{Odd characteristic}
We now consider the case when $\Fld_q$ has odd characteristic.   We first need a
well-known result.

\begin{lemma}
  If $\Fld_q$ has odd characteristic and $\nu\in \Fld_q$,
  then there exists $c,d\in\Fld_q$ such that $c^2+d^2=\nu$.
\end{lemma}

\begin{proof}
  Let $A=\{c^2 \st c\in \Fld_q\}$ and $B=\{\nu - d^2 \st d\in \Fld_q\}$.
  Since the map $\sigma\colon \Fld_q^\times \to \Fld_q^\times$ given by
  $\sigma\colon x\mapsto x^2$ has kernel $\{1,-1\}$, there are $(q-1)/2$
  squares in $\Fld_q\setminus \{0\}$.  Including zero, there are then
  $(q+1)/2$ squares in $\Fld_q$.  Thus $|A|=|B|=(q+1)/2$, so $A\cap
  B\neq \emptyset$, and $c^2=\nu-d^2$ for some $c,d\in\Fld_q$.
\end{proof}

Since there are $(q-1)/2$ nonzero squares in $\Fld_q$, given a nonsquare
$\nu\in\Fld_q$, the set $\{\nu b^2 \st b\in\Fld_q, b\neq 0\}$ is a set of
$(q-1)/2$ nonsquares in $\Fld_q$.  Consequently, every nonsquare is
equal to $\nu b^2$ for some $b\in\Fld_q$.

The matrix $aH$ for any $a\in \Fld_q$ is congruent to a diagonal matrix:
\begin{align*}
  \bmat 1&1\\-1&1\emat \bmat 0&a\\a&0\emat \bmat 1&-1\\1&1\emat
  = \bmat a&a\\a&-a\emat \bmat 1&-1\\1&1\emat = \bmat 2a&0\\0&-2a\emat.
\end{align*}
This fact combined with Lemma~\ref{lem:canonical-matrix} shows that
every symmetric matrix over $\Fld_q$ is congruent to a diagonal matrix.

\begin{lemma}
  Every invertible symmetric $k\times k$ matrix $B$ over $\Fld_q$ is
  congruent to either $I_k$ or $\diag(I_{k-1},\nu)$, where $\nu$ is
  any nonsquare in $\Fld_q$.
\end{lemma}

\begin{proof}
  Let $C$ be an invertible diagonal matrix congruent to $B$, with
  $C=N^tBN$, and let $\nu$ be any nonsquare in $\Fld_q$.

  By a permutation matrix $P$, let $D=P^tCP=\diag(b_1^2,b_2^2,\ldots,
  b_s^2, \nu c_1^2, \nu c_2^2,\ldots, \nu c_t^2)$, where the first $s$
  elements of the diagonal of $D$ are squares in $\Fld_q$ and the last
  $t$ elements are nonsquares in $\Fld_q$.

  Let
  $Q=\diag(b_1^{-1},b_2^{-1},\ldots,b_s^{-1},c_1^{-1},c_2^{-1},\ldots,c_t^{-1})$.
  Let $E=Q^tDQ=\diag(I_s,\nu I_t)$.

  Let $c,d\in \Fld_q$ such that $c^2+d^2=\nu$.  Let 
  \begin{align*}
R=\nu^{-1}\bmat
  c&d\\-d&c \emat.    
  \end{align*}
  Since $\det R=\nu^{-2}(c^2+d^2) = \nu^{-1} \neq 0$, $R$ is
  invertible.  Note that 
  \begin{align*}
    R^t(\nu I_2) R=\nu R^tR=\nu\nu^{-2}(c^2+d^2)I_2=I_2.
  \end{align*}
  If $t$ is even, let $S=\diag(I_s,R_1,R_2,\ldots,R_{t/2})$, where
  $R_i=R$ for each $i$.  Then $S^tES=I_k$.  If $t$ is odd, let
  $S=\diag(I_s,R_1,R_2,\ldots,R_{(t-1)/2},1)$.  Then
  $S^tES=\diag(I_{k-1},\nu)$.
\end{proof}

The next lemma shows that these two cases are in fact different and
gives a simple criteria to determine which congruence class any
symmetric matrix is in.

\begin{lemma}
  If $\det B$ is a square (nonsquare) and  $\hat B$ is congruent to
  $B$, then $\det \hat B$ is a square (nonsquare).
\end{lemma}

\begin{proof}
  Let $\hat B=C^tBC$.  Then $\det \hat B=(\det C)^2(\det B)$.  Thus
  $\det B$ is a square if and only if $\det \hat B$ is a square.
\end{proof}

Since $\det I_k=1$ is a square and $\det (\diag(I_{k-1},\nu))=\nu$ is a
nonsquare, we can determine if a matrix is congruent to $I_k$ or
congruent to $\diag(I_{k-1},\nu)$ by whether the determinant is a square
or not.

It appears then that $|\mathcal{B}|=2$.  However, we can do better in
one case since we only are concerned with whether an entry of $U^tBU$
is zero or nonzero and not with the actual value of the entry.

\begin{definition}
  Let $B$ and $\hat B$ be matrices.  If $\hat B=dC^tBC$ for some
  invertible matrix $C$ and some nonzero constant $d$, then $B$ and
  $\hat B$ are \emph{projectively congruent}.
\end{definition}

Since multiplying by a nonzero constant preserves the zero/nonzero
pattern in a matrix over a field, if $B$ and $\hat B$ are projectively
congruent, then $U^tBU$ and $U^t\hat BU$ give isomorphic graphs.

\begin{lemma}
  If $k$ is odd, then an invertible symmetric $k\times k$ matrix is
  projectively congruent to $I_k$.
\end{lemma}

\begin{proof}
  Let $k=2\ell-1$.  We can see that $\det (\nu\diag(I_{k-1},\nu)) =
  \nu^{2\ell-1}\nu=\nu^{2\ell}$ is a square.  Thus $\diag(I_{k-1},\nu)$ is
  projectively congruent to $I_k$.
\end{proof}

The results in this subsection give us the following lemma.

\begin{lemma}
  \label{lem:B-for-odd-char}
  Let $q$ be odd.  To determine $\g_k(\Fld_q)$, we may take
  $\mathcal{B}$ as follows: if $k$ is odd, then $\mathcal{B}=\{I_k\}$;
  if $k$ is even, then $\mathcal{B}=\{I_k,\diag(I_{k-1},\nu)\}$, where
  $\nu$ is any nonsquare in $\Fld_q$
\end{lemma}

\subsection{Summary}
Combining Lemmas~\ref{lem:B-for-even-char} and
\ref{lem:B-for-odd-char}, the results of this section can be
summarized as the following theorem.
\begin{theorem}\label{thm:main-theorem}
  The set $\g_k(\Fld_q)$ is the set of graphs of the matrices in
  $\{U^tBU \st B\in\mathcal{B}\}$, where the columns of $U$ are a
  maximal set of nonzero vectors in $\Fld_q^k$ such that no vector is
  a multiple of another and $\mathcal{B}$ is given by:
  \begin{enumerate}
  \item \label{case:k-odd}if $k$ is odd, $\mathcal{B}=\{I_k\}$.
  \item \label{case:k-even-q-even}if $k$ is even and $\characteristic \Fld_q=2$,
    $\mathcal{B}=\{I_k,\diag(H_1,H_2,\dots,H_{k/2})\}$, where
    $H_i=\left[\begin{matrix} 0&1\\1&0 \end{matrix}\right]$.
  \item \label{case:k-even-q-odd}if $k$ is even and $\characteristic \Fld_q\neq 2$,
    $\mathcal{B}=\{I_k,\diag(I_{k-1},\nu)\}$, where $\nu$ is any
    non-square in $\Fld_q$.
  \end{enumerate}
\end{theorem}

\subsection{Examples of characterizations}
\label{sec:example-subs-theorems}
As special cases of Theorem~\ref{thm:main-theorem}, we present the
following corollaries which calculate $\g_k(\Fld_q)$ for several
$k$ and $q$.  In the corollaries, we label a graph in
$\g_k(\Fld_q)$ using the pattern $FqRk$, signifying that it is a graph
for the $\mr(\Fld_q,G)\leq k$ corollary.  To compute these graphs, we
used the software program Sage \cite{sage-2.8.15} and the Sage functions
listed in Appendix~\ref{app:gen-subs-graphs}.

In these theorems, recall that $K_1$ does not have a loop.

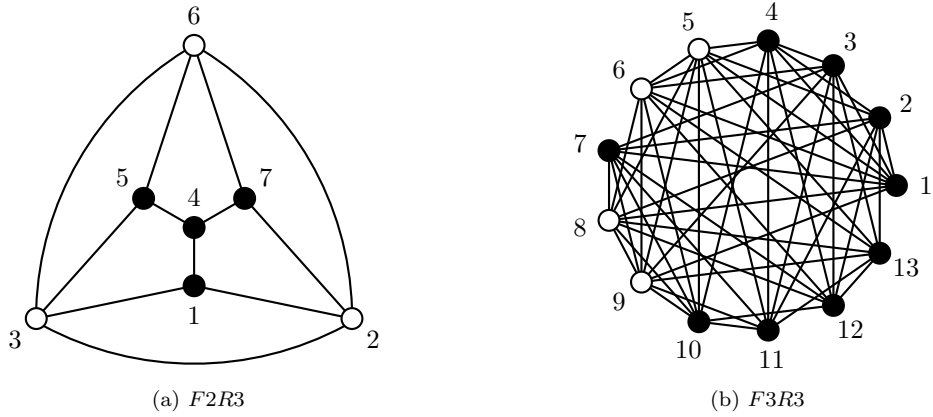
\begin{figure}
  \centering
 \subfloat[$F2R3$]{\begin{tikzpicture}[scale=2.52]
\ifx\dpiclw\undefined\newdimen\dpiclw\fi
\global\def\dpicdraw{\draw[line width=\dpiclw]}
\global\def\dpicstop{;}
\dpiclw=0.8bp
\dpicdraw[fill=black](0.943945,-0.199934) circle (0.021654in)\dpicstop
\dpicdraw[fill=black](0.677475,-0.046088) circle (0.021654in)\dpicstop
\dpicdraw (0.111228,-0.680703) circle (0.021654in)\dpicstop
\dpicdraw[fill=black](0.943945,-0.507627) circle (0.021654in)\dpicstop
\dpicdraw (1.776661,-0.680704) circle (0.021654in)\dpicstop
\dpicdraw[fill=black](1.210414,-0.046088) circle (0.021654in)\dpicstop
\dpicdraw (0.943945,0.761604) circle (0.021654in)\dpicstop
\draw (0.943945,-0.042634) node{4};
\draw (0.566247,0.06514) node{5};
\draw (-0,-0.791931) node{3};
\draw (0.943945,-0.664927) node{1};
\draw (1.887889,-0.791932) node{2};
\draw (1.321642,0.06514) node{7};
\draw (0.943945,0.918904) node{6};
\dpicdraw (0.896313,-0.172434)
 --(0.725107,-0.073588)\dpicstop
\dpicdraw (0.943945,-0.254934)
 --(0.943945,-0.452627)\dpicstop
\dpicdraw (0.991576,-0.172434)
 --(1.162783,-0.073588)\dpicstop
\dpicdraw (0.640858,-0.087127)
 --(0.147845,-0.639665)\dpicstop
\dpicdraw (0.165077,-0.669511)
 --(0.890096,-0.518819)\dpicstop
\dpicdraw (0.997794,-0.518819)
 --(1.722812,-0.669511)\dpicstop
\dpicdraw (1.740044,-0.639665)
 --(1.247032,-0.087127)\dpicstop
\dpicdraw (1.193182,0.006143)
 --(0.961177,0.709373)\dpicstop
\dpicdraw (0.926713,0.709373)
 --(0.694707,0.006143)\dpicstop
\dpicdraw (0.150119,-0.719594)
 ..controls (0.645975,-0.985341) and (1.241914,-0.985341)
 ..(1.73777,-0.719595)\dpicstop
\dpicdraw (0.990374,0.734851)
 ..controls (1.463283,0.441146) and (1.758287,-0.069316)
 ..(1.776661,-0.625704)\dpicstop
\dpicdraw (0.897516,0.734852)
 ..controls (0.424607,0.441147) and (0.129602,-0.069315)
 ..(0.111228,-0.625703)\dpicstop
\end{tikzpicture}\label{fig:f2r3}}
 \qquad\qquad\qquad
 \subfloat[$F3R3$]{\begin{tikzpicture}[scale=2.52]
\ifx\dpiclw\undefined\newdimen\dpiclw\fi
\global\def\dpicdraw{\draw[line width=\dpiclw]}
\global\def\dpicstop{;}
\dpiclw=0.8bp
\dpicdraw[fill=black](1.668838,0) circle (0.021654in)\dpicstop
\dpicdraw[fill=black](1.580727,0.357479) circle (0.021654in)\dpicstop
\dpicdraw[fill=black](1.33658,0.633065) circle (0.021654in)\dpicstop
\dpicdraw[fill=black](0.992328,0.763622) circle (0.021654in)\dpicstop
\dpicdraw (0.626834,0.719243) circle (0.021654in)\dpicstop
\dpicdraw (0.32383,0.510094) circle (0.021654in)\dpicstop
\dpicdraw[fill=black](0.152729,0.184089) circle (0.021654in)\dpicstop
\dpicdraw (0.152729,-0.184089) circle (0.021654in)\dpicstop
\dpicdraw (0.32383,-0.510094) circle (0.021654in)\dpicstop
\dpicdraw[fill=black](0.626834,-0.719243) circle (0.021654in)\dpicstop
\dpicdraw[fill=black](0.992328,-0.763622) circle (0.021654in)\dpicstop
\dpicdraw[fill=black](1.33658,-0.633065) circle (0.021654in)\dpicstop
\dpicdraw[fill=black](1.580727,-0.357479) circle (0.021654in)\dpicstop
\draw (1.826138,0) node{1};
\draw (1.72001,0.43058) node{2};
\draw (1.425937,0.76252) node{3};
\draw (1.011288,0.919775) node{4};
\draw (0.571055,0.866321) node{5};
\draw (0.206089,0.614404) node{6};
\draw (0,0.221733) node{7};
\draw (0,-0.221733) node{8};
\draw (0.206089,-0.614404) node{9};
\draw (0.571055,-0.866321) node{10};
\draw (1.011288,-0.919775) node{11};
\draw (1.425937,-0.76252) node{12};
\draw (1.72001,-0.43058) node{13};
\dpicdraw (1.655676,0.053402)
 --(1.59389,0.304078)\dpicstop
\dpicdraw (1.643278,0.0487)
 --(1.36214,0.584364)\dpicstop
\dpicdraw (1.632366,0.041168)
 --(1.0288,0.722454)\dpicstop
\dpicdraw (1.623574,0.031244)
 --(0.672099,0.688)\dpicstop
\dpicdraw (1.617412,0.019503)
 --(0.375256,0.490591)\dpicstop
\dpicdraw (1.614239,0.00663)
 --(0.207328,0.177459)\dpicstop
\dpicdraw (1.614239,-0.00663)
 --(0.207328,-0.177459)\dpicstop
\dpicdraw (1.617412,-0.019503)
 --(0.375256,-0.490591)\dpicstop
\dpicdraw (1.535463,0.388723)
 --(1.037592,0.732379)\dpicstop
\dpicdraw (1.529302,0.376983)
 --(0.67826,0.69974)\dpicstop
\dpicdraw (1.526129,0.35085)
 --(0.207328,0.190718)\dpicstop
\dpicdraw (1.529302,0.337976)
 --(0.204155,-0.164586)\dpicstop
\dpicdraw (1.555168,0.308779)
 --(1.017888,-0.714922)\dpicstop
\dpicdraw (1.567565,0.304078)
 --(1.349743,-0.579663)\dpicstop
\dpicdraw (1.580727,0.302479)
 --(1.580727,-0.302479)\dpicstop
\dpicdraw (1.285154,0.652568)
 --(1.043754,0.744119)\dpicstop
\dpicdraw (1.281981,0.626435)
 --(0.378429,0.516724)\dpicstop
\dpicdraw (1.285154,0.613561)
 --(0.204155,0.203592)\dpicstop
\dpicdraw (1.300109,0.591896)
 --(0.360302,-0.468926)\dpicstop
\dpicdraw (1.323418,0.579663)
 --(1.00549,-0.71022)\dpicstop
\dpicdraw (1.33658,0.578065)
 --(1.33658,-0.578065)\dpicstop
\dpicdraw (1.349743,0.579663)
 --(1.567565,-0.304078)\dpicstop
\dpicdraw (0.937729,0.756993)
 --(0.681433,0.725873)\dpicstop
\dpicdraw (0.940902,0.744119)
 --(0.375256,0.529598)\dpicstop
\dpicdraw (0.979166,0.71022)
 --(0.639997,-0.665841)\dpicstop
\dpicdraw (0.992328,0.708622)
 --(0.992328,-0.708622)\dpicstop
\dpicdraw (1.00549,0.71022)
 --(1.323418,-0.579663)\dpicstop
\dpicdraw (0.58157,0.688)
 --(0.369094,0.541338)\dpicstop
\dpicdraw (0.601275,0.670543)
 --(0.178289,-0.135389)\dpicstop
\dpicdraw (0.613672,0.665841)
 --(0.336992,-0.456693)\dpicstop
\dpicdraw (0.626834,0.664243)
 --(0.626834,-0.664243)\dpicstop
\dpicdraw (0.639997,0.665841)
 --(0.979166,-0.71022)\dpicstop
\dpicdraw (0.663306,0.678075)
 --(1.544256,-0.316311)\dpicstop
\dpicdraw (0.310668,0.456693)
 --(0.165892,-0.130687)\dpicstop
\dpicdraw (0.32383,0.455094)
 --(0.32383,-0.455094)\dpicstop
\dpicdraw (0.336992,0.456693)
 --(0.613672,-0.665841)\dpicstop
\dpicdraw (0.360302,0.468926)
 --(1.300109,-0.591896)\dpicstop
\dpicdraw (0.369094,0.478851)
 --(1.535463,-0.326236)\dpicstop
\dpicdraw (0.152729,0.129089)
 --(0.152729,-0.129089)\dpicstop
\dpicdraw (0.165892,0.130687)
 --(0.310668,-0.456693)\dpicstop
\dpicdraw (0.178289,0.135389)
 --(0.601275,-0.670543)\dpicstop
\dpicdraw (0.189201,0.142921)
 --(0.955856,-0.722454)\dpicstop
\dpicdraw (0.197993,0.152845)
 --(1.291316,-0.601821)\dpicstop
\dpicdraw (0.178289,-0.232789)
 --(0.29827,-0.461394)\dpicstop
\dpicdraw (0.189201,-0.225257)
 --(0.590363,-0.678075)\dpicstop
\dpicdraw (0.204155,-0.203592)
 --(1.285154,-0.613561)\dpicstop
\dpicdraw (0.207328,-0.190718)
 --(1.526129,-0.35085)\dpicstop
\dpicdraw (0.369094,-0.541338)
 --(0.58157,-0.688)\dpicstop
\dpicdraw (0.375256,-0.529598)
 --(0.940902,-0.744119)\dpicstop
\dpicdraw (0.378429,-0.503465)
 --(1.526129,-0.364109)\dpicstop
\dpicdraw (0.681433,-0.725873)
 --(0.937729,-0.756993)\dpicstop
\dpicdraw (0.681433,-0.712614)
 --(1.281981,-0.639694)\dpicstop
\dpicdraw (1.037592,-0.732379)
 --(1.535463,-0.388723)\dpicstop
\dpicdraw (1.373052,-0.591896)
 --(1.544256,-0.398647)\dpicstop
\end{tikzpicture}\label{fig:f3r3}}

 \caption{Graphs in Corollaries~\ref{cor:f2r3} and \ref{cor:f3r3}}
 \label{fig:r3}
\end{figure}

\begin{corollary}\label{cor:f2r3}
  Let $G$ be any simple graph.  Let $F2R3$ be the graph in
  Figure~\ref{fig:r3}\subref{fig:f2r3}.  Then $\mr(\Fld_2,G) \leq 3$ (i.e., $G\in
  \G_3(\Fld_2)$) if and only if $G$ is a blowup graph of $F2R3\cup
  K_1$.
\end{corollary}

\begin{proof}
  As matrices over $\Fld_2$, let
  \begin{equation*}
    U=\left[\begin{array}{rrrrrrr}
0 & 1 & 0 & 1 & 0 & 1 & 1 \\
0 & 0 & 1 & 1 & 1 & 1 & 0 \\
1 & 1 & 1 & 1 & 0 & 0 & 0
\end{array}\right]
 \quad \text{ and } \quad 
B=\left[\begin{array}{rrr}
1 & 0 & 0 \\
0 & 1 & 0 \\
0 & 0 & 1
\end{array}\right].
  \end{equation*}
  Then the graph $F2R3$ corresponds to the matrix
  \begin{equation*}
    U^tBU=          	
\left[\begin{array}{rrrrrrr}
1 & 1 & 1 & 1 & 0 & 0 & 0 \\
1 & 0 & 1 & 0 & 0 & 1 & 1 \\
1 & 1 & 0 & 0 & 1 & 1 & 0 \\
1 & 0 & 0 & 1 & 1 & 0 & 1 \\
0 & 0 & 1 & 1 & 1 & 1 & 0 \\
0 & 1 & 1 & 0 & 1 & 0 & 1 \\
0 & 1 & 0 & 1 & 0 & 1 & 1
\end{array}\right].
  \end{equation*}
\end{proof}

\begin{corollary} \label{cor:f3r3} Let $G$ be any simple graph.  Let
  $F3R3$ be the graph in Figure~\ref{fig:r3}\subref{fig:f3r3}. Then
  $\mr(\Fld_3,G) \leq 3$ (i.e., $G\in \G_3(\Fld_3)$) if and only if $G$ is
  a blowup graph of $F3R3\cup K_1$.
\end{corollary}

\begin{proof}
  As matrices over $\Fld_3$, let
  \begin{equation*}
    U=\left[\begin{array}{rrrrrrrrrrrrr}
0 & 1 & 2 & 0 & 1 & 2 & 0 & 1 & 2 & 0 & 1 & 2 & 1 \\
0 & 0 & 0 & 1 & 1 & 1 & 2 & 2 & 2 & 1 & 1 & 1 & 0 \\
1 & 1 & 1 & 1 & 1 & 1 & 1 & 1 & 1 & 0 & 0 & 0 & 0
\end{array}\right]
 \quad \text{ and } \quad 
B=\left[\begin{array}{rrr}
1 & 0 & 0 \\
0 & 1 & 0 \\
0 & 0 & 1
\end{array}\right].
  \end{equation*}
  Then the graph $F3R3$ corresponds to the matrix
  \begin{equation*}
    U^tBU= \left[\begin{array}{rrrrrrrrrrrrr}
1 & 1 & 1 & 1 & 1 & 1 & 1 & 1 & 1 & 0 & 0 & 0 & 0 \\
1 & 2 & 0 & 1 & 2 & 0 & 1 & 2 & 0 & 0 & 1 & 2 & 1 \\
1 & 0 & 2 & 1 & 0 & 2 & 1 & 0 & 2 & 0 & 2 & 1 & 2 \\
1 & 1 & 1 & 2 & 2 & 2 & 0 & 0 & 0 & 1 & 1 & 1 & 0 \\
1 & 2 & 0 & 2 & 0 & 1 & 0 & 1 & 2 & 1 & 2 & 0 & 1 \\
1 & 0 & 2 & 2 & 1 & 0 & 0 & 2 & 1 & 1 & 0 & 2 & 2 \\
1 & 1 & 1 & 0 & 0 & 0 & 2 & 2 & 2 & 2 & 2 & 2 & 0 \\
1 & 2 & 0 & 0 & 1 & 2 & 2 & 0 & 1 & 2 & 0 & 1 & 1 \\
1 & 0 & 2 & 0 & 2 & 1 & 2 & 1 & 0 & 2 & 1 & 0 & 2 \\
0 & 0 & 0 & 1 & 1 & 1 & 2 & 2 & 2 & 1 & 1 & 1 & 0 \\
0 & 1 & 2 & 1 & 2 & 0 & 2 & 0 & 1 & 1 & 2 & 0 & 1 \\
0 & 2 & 1 & 1 & 0 & 2 & 2 & 1 & 0 & 1 & 0 & 2 & 2 \\
0 & 1 & 2 & 0 & 1 & 2 & 0 & 1 & 2 & 0 & 1 & 2 & 1
\end{array}\right].
  \end{equation*}
\end{proof}

The next corollary gives the simplest previously-unknown result for
which $\g_k(\Fld_q)$ contains two graphs.

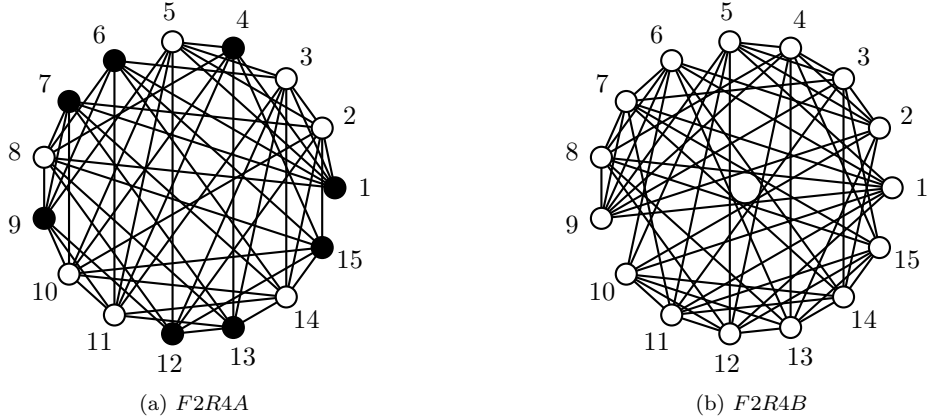
\begin{figure}[h!]
  \centering
 \subfloat[$F2R4A$]{\begin{tikzpicture}[scale=2.54]
\ifx\dpiclw\undefined\newdimen\dpiclw\fi
\global\def\dpicdraw{\draw[line width=\dpiclw]}
\global\def\dpicstop{;}
\dpiclw=0.8bp
\dpicdraw[fill=black](1.675515,0) circle (0.021654in)\dpicstop
\dpicdraw (1.609011,0.312874) circle (0.021654in)\dpicstop
\dpicdraw (1.421,0.57165) circle (0.021654in)\dpicstop
\dpicdraw[fill=black](1.143989,0.731582) circle (0.021654in)\dpicstop
\dpicdraw (0.825877,0.765017) circle (0.021654in)\dpicstop
\dpicdraw[fill=black](0.521668,0.666173) circle (0.021654in)\dpicstop
\dpicdraw[fill=black](0.283963,0.452143) circle (0.021654in)\dpicstop
\dpicdraw (0.153863,0.159932) circle (0.021654in)\dpicstop
\dpicdraw[fill=black](0.153863,-0.159932) circle (0.021654in)\dpicstop
\dpicdraw (0.283963,-0.452143) circle (0.021654in)\dpicstop
\dpicdraw (0.521668,-0.666173) circle (0.021654in)\dpicstop
\dpicdraw[fill=black](0.825877,-0.765017) circle (0.021654in)\dpicstop
\dpicdraw[fill=black](1.143989,-0.731582) circle (0.021654in)\dpicstop
\dpicdraw (1.421,-0.57165) circle (0.021654in)\dpicstop
\dpicdraw[fill=black](1.609011,-0.312874) circle (0.021654in)\dpicstop
\draw (1.832815,0) node{1};
\draw (1.752712,0.376854) node{2};
\draw (1.526254,0.688547) node{3};
\draw (1.192598,0.881183) node{4};
\draw (0.809435,0.921455) node{5};
\draw (0.443018,0.802399) node{6};
\draw (0.156705,0.544601) node{7};
\draw (0,0.192637) node{8};
\draw (0,-0.192637) node{9};
\draw (0.156705,-0.544601) node{10};
\draw (0.443018,-0.802399) node{11};
\draw (0.809435,-0.921455) node{12};
\draw (1.192598,-0.881183) node{13};
\draw (1.526254,-0.688547) node{14};
\draw (1.752712,-0.376854) node{15};
\dpicdraw (1.664079,0.053798)
 --(1.620446,0.259076)\dpicstop
\dpicdraw (1.653144,0.050245)
 --(1.44337,0.521405)\dpicstop
\dpicdraw (1.643186,0.044496)
 --(1.176317,0.687086)\dpicstop
\dpicdraw (1.634642,0.036802)
 --(0.86675,0.728215)\dpicstop
\dpicdraw (1.627883,0.0275)
 --(0.5693,0.638673)\dpicstop
\dpicdraw (1.623207,0.016996)
 --(0.336271,0.435147)\dpicstop
\dpicdraw (1.620816,0.005749)
 --(0.208561,0.154183)\dpicstop
\dpicdraw (1.576683,0.35737)
 --(1.453328,0.527154)\dpicstop
\dpicdraw (1.56138,0.340374)
 --(0.873509,0.737517)\dpicstop
\dpicdraw (1.554312,0.318623)
 --(0.338662,0.446393)\dpicstop
\dpicdraw (1.56138,0.285374)
 --(0.331594,-0.424643)\dpicstop
\dpicdraw (1.576683,0.268378)
 --(0.858206,-0.720521)\dpicstop
\dpicdraw (1.597576,0.259076)
 --(1.432435,-0.517852)\dpicstop
\dpicdraw (1.609011,0.257874)
 --(1.609011,-0.257874)\dpicstop
\dpicdraw (1.368692,0.588646)
 --(0.878185,0.748021)\dpicstop
\dpicdraw (1.366301,0.577399)
 --(0.576367,0.660424)\dpicstop
\dpicdraw (1.388672,0.527154)
 --(0.553997,-0.621677)\dpicstop
\dpicdraw (1.398629,0.521405)
 --(0.848248,-0.714772)\dpicstop
\dpicdraw (1.409565,0.517852)
 --(1.155424,-0.677784)\dpicstop
\dpicdraw (1.421,0.51665)
 --(1.421,-0.51665)\dpicstop
\dpicdraw (1.089291,0.737331)
 --(0.880576,0.759268)\dpicstop
\dpicdraw (1.096358,0.704082)
 --(0.201494,0.187432)\dpicstop
\dpicdraw (1.111661,0.687086)
 --(0.316291,-0.407647)\dpicstop
\dpicdraw (1.121619,0.681337)
 --(0.544039,-0.615928)\dpicstop
\dpicdraw (1.143989,0.676582)
 --(1.143989,-0.676582)\dpicstop
\dpicdraw (1.16636,0.681337)
 --(1.586641,-0.262629)\dpicstop
\dpicdraw (0.793549,0.720521)
 --(0.186191,-0.115436)\dpicstop
\dpicdraw (0.803507,0.714772)
 --(0.306334,-0.401898)\dpicstop
\dpicdraw (0.814442,0.711219)
 --(0.533104,-0.612375)\dpicstop
\dpicdraw (0.825877,0.710017)
 --(0.825877,-0.710017)\dpicstop
\dpicdraw (0.48934,0.621677)
 --(0.186191,0.204428)\dpicstop
\dpicdraw (0.499298,0.615928)
 --(0.176233,-0.109687)\dpicstop
\dpicdraw (0.521668,0.611173)
 --(0.521668,-0.611173)\dpicstop
\dpicdraw (0.553997,0.621677)
 --(1.388672,-0.527154)\dpicstop
\dpicdraw (0.562541,0.629371)
 --(1.568138,-0.276072)\dpicstop
\dpicdraw (0.261593,0.401898)
 --(0.176233,0.210177)\dpicstop
\dpicdraw (0.272528,0.398344)
 --(0.165298,-0.106134)\dpicstop
\dpicdraw (0.283963,0.397143)
 --(0.283963,-0.397143)\dpicstop
\dpicdraw (0.316291,0.407647)
 --(1.111661,-0.687086)\dpicstop
\dpicdraw (0.324836,0.41534)
 --(1.380127,-0.534848)\dpicstop
\dpicdraw (0.153863,0.104932)
 --(0.153863,-0.104932)\dpicstop
\dpicdraw (0.186191,0.115436)
 --(0.793549,-0.720521)\dpicstop
\dpicdraw (0.194736,0.12313)
 --(1.103116,-0.69478)\dpicstop
\dpicdraw (0.206171,0.142936)
 --(1.556703,-0.295878)\dpicstop
\dpicdraw (0.176233,-0.210177)
 --(0.261593,-0.401898)\dpicstop
\dpicdraw (0.186191,-0.204428)
 --(0.48934,-0.621677)\dpicstop
\dpicdraw (0.194736,-0.196734)
 --(0.785004,-0.728215)\dpicstop
\dpicdraw (0.324836,-0.488945)
 --(0.480795,-0.629371)\dpicstop
\dpicdraw (0.338662,-0.457892)
 --(1.366301,-0.565901)\dpicstop
\dpicdraw (0.338662,-0.446393)
 --(1.554312,-0.318623)\dpicstop
\dpicdraw (0.576367,-0.671922)
 --(1.089291,-0.725833)\dpicstop
\dpicdraw (0.576367,-0.660424)
 --(1.366301,-0.577399)\dpicstop
\dpicdraw (0.880576,-0.759268)
 --(1.089291,-0.737331)\dpicstop
\dpicdraw (0.873509,-0.737517)
 --(1.56138,-0.340374)\dpicstop
\dpicdraw (1.191621,-0.704082)
 --(1.373368,-0.59915)\dpicstop
\dpicdraw (1.453328,-0.527154)
 --(1.576683,-0.35737)\dpicstop
\end{tikzpicture}}
 \qquad\qquad\qquad
 \subfloat[$F2R4B$]{\begin{tikzpicture}[scale=2.54]
\ifx\dpiclw\undefined\newdimen\dpiclw\fi
\global\def\dpicdraw{\draw[line width=\dpiclw]}
\global\def\dpicstop{;}
\dpiclw=0.8bp
\dpicdraw (1.675515,0) circle (0.021654in)\dpicstop
\dpicdraw (1.609011,0.312874) circle (0.021654in)\dpicstop
\dpicdraw (1.421,0.57165) circle (0.021654in)\dpicstop
\dpicdraw (1.143989,0.731582) circle (0.021654in)\dpicstop
\dpicdraw (0.825877,0.765017) circle (0.021654in)\dpicstop
\dpicdraw (0.521668,0.666173) circle (0.021654in)\dpicstop
\dpicdraw (0.283963,0.452143) circle (0.021654in)\dpicstop
\dpicdraw (0.153863,0.159932) circle (0.021654in)\dpicstop
\dpicdraw (0.153863,-0.159932) circle (0.021654in)\dpicstop
\dpicdraw (0.283963,-0.452143) circle (0.021654in)\dpicstop
\dpicdraw (0.521668,-0.666173) circle (0.021654in)\dpicstop
\dpicdraw (0.825877,-0.765017) circle (0.021654in)\dpicstop
\dpicdraw (1.143989,-0.731582) circle (0.021654in)\dpicstop
\dpicdraw (1.421,-0.57165) circle (0.021654in)\dpicstop
\dpicdraw (1.609011,-0.312874) circle (0.021654in)\dpicstop
\draw (1.832815,0) node{1};
\draw (1.752712,0.376854) node{2};
\draw (1.526254,0.688547) node{3};
\draw (1.192598,0.881183) node{4};
\draw (0.809435,0.921455) node{5};
\draw (0.443018,0.802399) node{6};
\draw (0.156705,0.544601) node{7};
\draw (0,0.192637) node{8};
\draw (0,-0.192637) node{9};
\draw (0.156705,-0.544601) node{10};
\draw (0.443018,-0.802399) node{11};
\draw (0.809435,-0.921455) node{12};
\draw (1.192598,-0.881183) node{13};
\draw (1.526254,-0.688547) node{14};
\draw (1.752712,-0.376854) node{15};
\dpicdraw (1.634642,0.036802)
 --(0.86675,0.728215)\dpicstop
\dpicdraw (1.627883,0.0275)
 --(0.5693,0.638673)\dpicstop
\dpicdraw (1.623207,0.016996)
 --(0.336271,0.435147)\dpicstop
\dpicdraw (1.620816,0.005749)
 --(0.208561,0.154183)\dpicstop
\dpicdraw (1.620816,-0.005749)
 --(0.208561,-0.154183)\dpicstop
\dpicdraw (1.623207,-0.016996)
 --(0.336271,-0.435147)\dpicstop
\dpicdraw (1.627883,-0.0275)
 --(0.5693,-0.638673)\dpicstop
\dpicdraw (1.634642,-0.036802)
 --(0.86675,-0.728215)\dpicstop
\dpicdraw (1.576683,0.35737)
 --(1.453328,0.527154)\dpicstop
\dpicdraw (1.568138,0.349677)
 --(1.184862,0.69478)\dpicstop
\dpicdraw (1.56138,0.340374)
 --(0.873509,0.737517)\dpicstop
\dpicdraw (1.556703,0.32987)
 --(0.573977,0.649177)\dpicstop
\dpicdraw (1.556703,0.295878)
 --(0.206171,-0.142936)\dpicstop
\dpicdraw (1.56138,0.285374)
 --(0.331594,-0.424643)\dpicstop
\dpicdraw (1.586641,0.262629)
 --(1.16636,-0.681337)\dpicstop
\dpicdraw (1.597576,0.259076)
 --(1.432435,-0.517852)\dpicstop
\dpicdraw (1.373368,0.59915)
 --(1.191621,0.704082)\dpicstop
\dpicdraw (1.368692,0.588646)
 --(0.878185,0.748021)\dpicstop
\dpicdraw (1.366301,0.565901)
 --(0.338662,0.457892)\dpicstop
\dpicdraw (1.373368,0.54415)
 --(0.201494,-0.132432)\dpicstop
\dpicdraw (1.388672,0.527154)
 --(0.553997,-0.621677)\dpicstop
\dpicdraw (1.421,0.51665)
 --(1.421,-0.51665)\dpicstop
\dpicdraw (1.432435,0.517852)
 --(1.597576,-0.259076)\dpicstop
\dpicdraw (1.089291,0.737331)
 --(0.880576,0.759268)\dpicstop
\dpicdraw (1.096358,0.704082)
 --(0.201494,0.187432)\dpicstop
\dpicdraw (1.103116,0.69478)
 --(0.194736,-0.12313)\dpicstop
\dpicdraw (1.132554,0.677784)
 --(0.837312,-0.711219)\dpicstop
\dpicdraw (1.143989,0.676582)
 --(1.143989,-0.676582)\dpicstop
\dpicdraw (1.16636,0.681337)
 --(1.586641,-0.262629)\dpicstop
\dpicdraw (0.793549,0.720521)
 --(0.186191,-0.115436)\dpicstop
\dpicdraw (0.803507,0.714772)
 --(0.306334,-0.401898)\dpicstop
\dpicdraw (0.814442,0.711219)
 --(0.533104,-0.612375)\dpicstop
\dpicdraw (0.825877,0.710017)
 --(0.825877,-0.710017)\dpicstop
\dpicdraw (0.480795,0.629371)
 --(0.324836,0.488945)\dpicstop
\dpicdraw (0.48934,0.621677)
 --(0.186191,0.204428)\dpicstop
\dpicdraw (0.499298,0.615928)
 --(0.176233,-0.109687)\dpicstop
\dpicdraw (0.510233,0.612375)
 --(0.295398,-0.398344)\dpicstop
\dpicdraw (0.544039,0.615928)
 --(1.121619,-0.681337)\dpicstop
\dpicdraw (0.553997,0.621677)
 --(1.388672,-0.527154)\dpicstop
\dpicdraw (0.261593,0.401898)
 --(0.176233,0.210177)\dpicstop
\dpicdraw (0.272528,0.398344)
 --(0.165298,-0.106134)\dpicstop
\dpicdraw (0.295398,0.398344)
 --(0.510233,-0.612375)\dpicstop
\dpicdraw (0.324836,0.41534)
 --(1.380127,-0.534848)\dpicstop
\dpicdraw (0.331594,0.424643)
 --(1.56138,-0.285374)\dpicstop
\dpicdraw (0.153863,0.104932)
 --(0.153863,-0.104932)\dpicstop
\dpicdraw (0.186191,0.115436)
 --(0.793549,-0.720521)\dpicstop
\dpicdraw (0.194736,0.12313)
 --(1.103116,-0.69478)\dpicstop
\dpicdraw (0.206171,0.142936)
 --(1.556703,-0.295878)\dpicstop
\dpicdraw (0.324836,-0.488945)
 --(0.480795,-0.629371)\dpicstop
\dpicdraw (0.331594,-0.479643)
 --(0.778246,-0.737517)\dpicstop
\dpicdraw (0.336271,-0.469138)
 --(1.091681,-0.714586)\dpicstop
\dpicdraw (0.338662,-0.457892)
 --(1.366301,-0.565901)\dpicstop
\dpicdraw (0.573977,-0.683169)
 --(0.773569,-0.748021)\dpicstop
\dpicdraw (0.576367,-0.660424)
 --(1.366301,-0.577399)\dpicstop
\dpicdraw (0.573977,-0.649177)
 --(1.556703,-0.32987)\dpicstop
\dpicdraw (0.880576,-0.759268)
 --(1.089291,-0.737331)\dpicstop
\dpicdraw (0.873509,-0.737517)
 --(1.56138,-0.340374)\dpicstop
\dpicdraw (1.191621,-0.704082)
 --(1.373368,-0.59915)\dpicstop
\dpicdraw (1.184862,-0.69478)
 --(1.568138,-0.349677)\dpicstop
\dpicdraw (1.453328,-0.527154)
 --(1.576683,-0.35737)\dpicstop
\end{tikzpicture}}

 \caption{Graphs in Corollary~\ref{cor:f2r4}}
 \label{fig:f2r4}
\end{figure}

\begin{corollary} \label{cor:f2r4} Let $G$ be any simple graph.  Let
  $F2R4A$ and $F2R4B$ be the graphs in Figure~\ref{fig:f2r4}.  Then
  $\mr(\Fld_2,G) \leq 4$ (i.e., $G\in \G_4(\Fld_2)$) if and only if $G$ is
  a blowup graph of either $F2R4A\cup K_1$ or $F2R4B\cup K_1$.
\end{corollary}

\begin{proof}
  As matrices over $\Fld_2$, let
  \begin{equation*}
    U= \left[\begin{array}{rrrrrrrrrrrrrrr}
0 & 1 & 0 & 1 & 0 & 1 & 0 & 1 & 0 & 1 & 0 & 1 & 0 & 1 & 1 \\
0 & 0 & 1 & 1 & 0 & 0 & 1 & 1 & 0 & 0 & 1 & 1 & 1 & 1 & 0 \\
0 & 0 & 0 & 0 & 1 & 1 & 1 & 1 & 1 & 1 & 1 & 1 & 0 & 0 & 0 \\
1 & 1 & 1 & 1 & 1 & 1 & 1 & 1 & 0 & 0 & 0 & 0 & 0 & 0 & 0
\end{array}\right]
\end{equation*}
and let 
\begin{equation*}
B_1=I_4\quad \text{ and } \quad B_2=\left[\begin{array}{rrrr}
0 & 1 & 0 & 0 \\
1 & 0 & 0 & 0 \\
0 & 0 & 0 & 1 \\
0 & 0 & 1 & 0
\end{array}\right].
  \end{equation*}
  Then the graph $F2R4A$ corresponds to the matrix
  \begin{equation*}
    U^tB_1U=          	
\left[\begin{array}{rrrrrrrrrrrrrrr}
1 & 1 & 1 & 1 & 1 & 1 & 1 & 1 & 0 & 0 & 0 & 0 & 0 & 0 & 0 \\
1 & 0 & 1 & 0 & 1 & 0 & 1 & 0 & 0 & 1 & 0 & 1 & 0 & 1 & 1 \\
1 & 1 & 0 & 0 & 1 & 1 & 0 & 0 & 0 & 0 & 1 & 1 & 1 & 1 & 0 \\
1 & 0 & 0 & 1 & 1 & 0 & 0 & 1 & 0 & 1 & 1 & 0 & 1 & 0 & 1 \\
1 & 1 & 1 & 1 & 0 & 0 & 0 & 0 & 1 & 1 & 1 & 1 & 0 & 0 & 0 \\
1 & 0 & 1 & 0 & 0 & 1 & 0 & 1 & 1 & 0 & 1 & 0 & 0 & 1 & 1 \\
1 & 1 & 0 & 0 & 0 & 0 & 1 & 1 & 1 & 1 & 0 & 0 & 1 & 1 & 0 \\
1 & 0 & 0 & 1 & 0 & 1 & 1 & 0 & 1 & 0 & 0 & 1 & 1 & 0 & 1 \\
0 & 0 & 0 & 0 & 1 & 1 & 1 & 1 & 1 & 1 & 1 & 1 & 0 & 0 & 0 \\
0 & 1 & 0 & 1 & 1 & 0 & 1 & 0 & 1 & 0 & 1 & 0 & 0 & 1 & 1 \\
0 & 0 & 1 & 1 & 1 & 1 & 0 & 0 & 1 & 1 & 0 & 0 & 1 & 1 & 0 \\
0 & 1 & 1 & 0 & 1 & 0 & 0 & 1 & 1 & 0 & 0 & 1 & 1 & 0 & 1 \\
0 & 0 & 1 & 1 & 0 & 0 & 1 & 1 & 0 & 0 & 1 & 1 & 1 & 1 & 0 \\
0 & 1 & 1 & 0 & 0 & 1 & 1 & 0 & 0 & 1 & 1 & 0 & 1 & 0 & 1 \\
0 & 1 & 0 & 1 & 0 & 1 & 0 & 1 & 0 & 1 & 0 & 1 & 0 & 1 & 1
\end{array}\right]
\end{equation*}
and the graph $F2R4B$ corresponds to the matrix
\begin{equation*}
  U^tB_2U=\left[\begin{array}{rrrrrrrrrrrrrrr}
0 & 0 & 0 & 0 & 1 & 1 & 1 & 1 & 1 & 1 & 1 & 1 & 0 & 0 & 0 \\
0 & 0 & 1 & 1 & 1 & 1 & 0 & 0 & 1 & 1 & 0 & 0 & 1 & 1 & 0 \\
0 & 1 & 0 & 1 & 1 & 0 & 1 & 0 & 1 & 0 & 1 & 0 & 0 & 1 & 1 \\
0 & 1 & 1 & 0 & 1 & 0 & 0 & 1 & 1 & 0 & 0 & 1 & 1 & 0 & 1 \\
1 & 1 & 1 & 1 & 0 & 0 & 0 & 0 & 1 & 1 & 1 & 1 & 0 & 0 & 0 \\
1 & 1 & 0 & 0 & 0 & 0 & 1 & 1 & 1 & 1 & 0 & 0 & 1 & 1 & 0 \\
1 & 0 & 1 & 0 & 0 & 1 & 0 & 1 & 1 & 0 & 1 & 0 & 0 & 1 & 1 \\
1 & 0 & 0 & 1 & 0 & 1 & 1 & 0 & 1 & 0 & 0 & 1 & 1 & 0 & 1 \\
1 & 1 & 1 & 1 & 1 & 1 & 1 & 1 & 0 & 0 & 0 & 0 & 0 & 0 & 0 \\
1 & 1 & 0 & 0 & 1 & 1 & 0 & 0 & 0 & 0 & 1 & 1 & 1 & 1 & 0 \\
1 & 0 & 1 & 0 & 1 & 0 & 1 & 0 & 0 & 1 & 0 & 1 & 0 & 1 & 1 \\
1 & 0 & 0 & 1 & 1 & 0 & 0 & 1 & 0 & 1 & 1 & 0 & 1 & 0 & 1 \\
0 & 1 & 0 & 1 & 0 & 1 & 0 & 1 & 0 & 1 & 0 & 1 & 0 & 1 & 1 \\
0 & 1 & 1 & 0 & 0 & 1 & 1 & 0 & 0 & 1 & 1 & 0 & 1 & 0 & 1 \\
0 & 0 & 1 & 1 & 0 & 0 & 1 & 1 & 0 & 0 & 1 & 1 & 1 & 1 & 0
\end{array}\right].
\end{equation*}
\end{proof}

\section{Connection to projective geometry}

As mentioned previously, the classifications of symmetric matrices in
Section~\ref{sec:congr-class-symm} are standard classification results
in projective geometry.  In this section, we first review
appropriate terminology and highlight this connection to projective
geometry.  We will define slightly more terminology than is strictly
necessary to help the reader see where these things fit into standard
projective geometry.  We then give some examples of how results in
projective geometry can help us understand $\g_k(\Fld_q)$ better.  For
further material, a definitive treatise on projective geometry is
contained in the series \cite{hirschfeld-projective} and
\cite{hirschfeld-general-geometries}.

\subsection{Definitions and the connection}
We start with basic definitions from projective geometry.

\begin{definition}
  Let $V=\Fld_q^{n+1}$, the vector space of dimension $n+1$ over
  $\Fld_q$.  For $x,y\in V-\vec 0 $, we define an equivalence
  relation by
\begin{align*}
x\sim y \iff x=cy,\quad \text{ where } c\in\Fld_q \text{ and } c\neq 0.
\end{align*}
Denote the equivalence class containing $x\in V-\0$ as $\x=\{cx \st
c\in \Fld_q \text{ and } c\neq 0 \}$.  Geometrically, we can think of
the class $\x$ as the set of non-origin points on a line passing
through $x$ and the origin in $V$.  These equivalence classes form the
\emph{projective geometry} $PG(n,q)$ of (projective) dimension $n$ and order $q$.
The equivalence classes are called the \emph{points} of $PG(n,q)$.
Each subspace of dimension $m+1$ in $V$ corresponds to a subspace of
(projective) dimension $m$ in $PG(n,q)$.  If a projective geometry has
(projective) dimension 2, then it is called a \emph{projective plane}.
\end{definition}

Note that there is a shift by one in dimension between a vector space
$V$ and its subspaces and the projective geometry associated with $V$
and its subspaces.  To help the reader, we will use the nonstandard term
\emph{projective dimension} (or ``$\projdim$'') when dealing with the
dimension of a projective geometry.

\begin{definition}
  Let $\mathcal{S}$ be the set of subspaces of
  $PG(n,q)$.  A \emph{correlation} $\sigma: \mathcal{S} \to
  \mathcal{S}$ is a bijective map such that for any subspaces
  $R,T\in\mathcal{S}$, $R\subseteq T$ implies that $\sigma(T)\subseteq
  \sigma(R)$ and $\projdim \sigma(R)=n-1-\projdim R$.  A
  \emph{polarity} is a correlation $\sigma$ of order 2 (i.e.,
  $\sigma^2=1$, the identity map).
\end{definition}

Note that any polarity $\sigma$ maps points in $\mathcal{S}$ to hyperplanes
(subspaces of projective dimension $n-1$ in $\mathcal{S}$) and hyperplanes to
points.  Since $\sigma^2=1$, we have $Y=\sigma(\x)$ if and only if $\sigma(Y)=\x$, so 
$\sigma$ induces a bijection between points and hyperplanes. This
bijection leads to the next definition.

\begin{definition}
  Let $\sigma$ be a polarity on $PG(n,q)$.  Let $\x,\y$ be points in
  $PG(n,q)$.  We say that $\sigma(\x)$ is the \emph{polar}
  (hyperplane) of $\x$ and $\x$ is the \emph{pole} of $\sigma(\x)$.
  If $\y\in\sigma(\x)$, then $\x\in\sigma(\y)$ and we say that $\x$
  and $\y$ are \emph{conjugate} points.  If $\x\in\sigma(\x)$, then we
  say that $\x$ is \emph{self-conjugate} or \emph{absolute}.
  Similarly, if $S$ is a subspace of $PG(n,q)$, then $S$ is
  \emph{absolute} if $\sigma(S)\subseteq S$ or $S\subseteq\sigma(S)$.
  A subspace of $PG(n,q)$ consisting of absolute points is called
  \emph{isotropic}.
\end{definition}

The next definition gives the connection with symmetric matrices.

\begin{definition}
  Let $B$ be an $(n+1)\times (n+1)$ invertible symmetric matrix over
  $\Fld_q$.  Define $\sigma:\mathcal{S} \to \mathcal{S}$ by $\sigma: R\mapsto R^\perp$,
  where the orthogonality relation is defined by the nondegenerate
  symmetric bilinear form represented by $B$ (i.e., $R^\perp=\{\y \st
  x^tBy=0 \text{ for all } \x\in R\}$).  We call $\sigma$ the polarity
  associated with $B$.
\end{definition}

The fact that the $\sigma$ in the previous definition is a polarity is
easy to check.

Let $M_1$ and $M_2$ be symmetric matrices. Let $\sigma_1$ and
$\sigma_2$ be the associated polarities, respectively.  Two polarities
are equivalent if the matrices are projectively congruent, i.e.,
$\sigma_1$ is equivalent to $\sigma_2$ if $M_1=dC^tM_2C$ for some
nonzero $d$ and invertible matrix $C$.  Thus there is a unique
polarity associated with each matrix given in 
Theorem~\ref{thm:main-theorem}.

We now summarize from \cite[Section 2.1.5]{hirschfeld-projective} the
classification of polarities that are associated with symmetric
matrices.  Let $B$ be an invertible symmetric matrix over $\Fld_q$.
Let $\sigma$ be the polarity associated with $B$.
\begin{itemize}
\item If $q$ is odd, then $\sigma$ is called an \emph{ordinary
    polarity}.  

  If $B$ has even order, then the associated polarity is either a 
  \emph{hyperbolic} polarity or an \emph{elliptic} polarity.  The
  correspondence between these types of polarities and the matrices in
  $\mathcal{B}$ from
  Theorem~\ref{thm:main-theorem}(\ref{case:k-even-q-odd}) is slightly
  nontrivial and is summarized in
  \cite[Corollary~5.19]{hirschfeld-projective}.

  If $B$ has odd order, then $\sigma$ is a \emph{parabolic} polarity, which
  corresponds to $\mathcal{B}$ in
  Theorem~\ref{thm:main-theorem}(\ref{case:k-odd}).

\item If $q$ is even and $b_{ii}=0$ for all $i$, then $\sigma$ is a
  \emph{null} polarity (or in alternate terminology, $\sigma$ is a
  \emph{symplectic} polarity).  Note that this only occurs when $B$
  has even order since otherwise $B$ is not invertible.  This case
  corresponds to the non-identity matrix in the $\mathcal{B}$ in
  Theorem~\ref{thm:main-theorem}(\ref{case:k-even-q-even}).

\item If $q$ is even and there is some $b_{ii}\neq 0$, then $\sigma$ is a
  \emph{pseudo-polarity}. This case corresponds to the identity matrix
  in $\mathcal{B}$ in Theorem~\ref{thm:main-theorem}(\ref{case:k-odd})
  or (\ref{case:k-even-q-even}).
\end{itemize}

We pause to note that there are polarities that are not associated
with symmetric matrices.  However, since we are only concerned about
symmetric matrices, we will restrict ourselves to this case.
Information about polarities not associated with symmetric matrices
may also be found in \cite{hirschfeld-projective}.

We now examine the connection to graphs by recalling the definition of
a polarity graph.

\begin{definition}
  Let $B$ be an invertible symmetric $(n+1)\times (n+1)$ matrix over $\Fld_q$ and
  let $\sigma$ be the associated polarity.  The \emph{polarity graph} of
  $\sigma$ has as its vertices the points of $PG(n,q)$ and as its edges
  $\{\x\y \st x^tBy=0\}$.
\end{definition}
In a polarity graph, $\x$ is adjacent to $\y$ exactly when $\x$ and
$\y$ are conjugate (i.e., $x$ and $y$ are orthogonal with respect to
$B$).  In standard literature, loops are not allowed in polarity
graphs.  However, for our purposes, loops convey needed information,
so a vertex $\x$ in a polarity graph has a loop if and only if $\x$ is
absolute (i.e., $x^tBx=0$, where $B$ is an invertible symmetric matrix
associated with the polarity).

In Theorem~\ref{thm:main-theorem}, the vertices of a graph in
$\g_k(\Fld_q)$ represent the points of the projective geometry
$PG(k-1,q)$ and an edge is drawn if the corresponding points are
\emph{not} conjugate (i.e., $x^tBy\neq 0$).  Thus, the graphs in
Theorem~\ref{thm:main-theorem} are exactly the complements of polarity
graphs.  Recall that when dealing with looped graphs, a vertex is
looped in the complement of a graph if and only if it is nonlooped in
the original graph.

Using this connection, we can restate Theorem~\ref{thm:main-theorem}:
\begin{theorem}\label{thm:main-theorem-proj}
  The set $\g_k(\Fld_q)$ is the set of complements of the (looped)
  polarity graphs of the polarities on $PG(k-1,q)$ that are associated
  with symmetric matrices.
\end{theorem}

\subsection{Consequences of the connection}

With the main theorem stated as in
Theorem~\ref{thm:main-theorem-proj}, we can use a variety of known
results about polarity graphs to derive results about graphs in
$\g_k(\Fld_q)$.  In this section, we list a few consequences of
Theorem~\ref{thm:main-theorem-proj}.

An elementary result in projective geometry gives us
the size of the graphs in $\g_k(\Fld_q)$.  While this result could
have been realized from the statement in
Theorem~\ref{thm:main-theorem}, it also naturally follows as a
consequence of Theorem~\ref{thm:main-theorem-proj}.
\begin{theorem}
  Every graph in $\g_k(\Fld_q)$ has $\frac{q^k-1}{q-1}$ vertices.
\end{theorem}

\begin{proof}
  There are $q^k-1$ vectors in $\Fld_q^k-\0$.  Since there are $q-1$
  nonzero constants in $\Fld_q$, there are $q-1$ elements in each
  equivalence class in $PG(k-1,q)$, so there are $\frac{q^k-1}{q-1}$
  points in $PG(k-1,q)$.
\end{proof}

The following observation follows directly from
Theorem~\ref{thm:main-theorem-proj} and restates the criteria for an
edge in a graph in $\g_k(\Fld_q)$ in several ways.

\begin{observation}\label{obs:edges}
  Let $G\in\g_k(\Fld_q)$ and let $u$ and $v$ be (not necessarily
  distinct) vertices in $G$.  Let $\sigma$ be the polarity
  corresponding to $G$ and let $B$ be an invertible symmetric matrix
  corresponding to $\sigma$. Then $uv$ is an edge in $G$ if and only
  if:
  \begin{enumerate}
  \item $u^tBv\neq 0$ (equivalently, $v^tBu\neq 0$), or equivalently,
  \item $u$ and $v$ are not conjugate points, or equivalently,
  \item $u\not\in \sigma(v)$ (equivalently, $v\not\in\sigma(u)$).
  \end{enumerate}
\end{observation}

\begin{corollary}\label{cor:regular}
  A graph $G\in\g_k(\Fld_q)$ is regular of degree $q^{k-1}$ (using the
  convention that a loop adds one to the degree of a vertex).
\end{corollary}
\begin{proof}
  Let $v\in G$ and let $\sigma$ be the polarity associated with $G$.  Since
  the hyperplane $\sigma(v)$ contains $\frac{q^{k-1}-1}{q-1}$ points, this
  is the degree of a $v$ in the complement of $G$.  Thus the degree of
  $v$ in $G$ is
  \begin{equation*}
    \frac{q^{k}-1}{q-1} - \frac{q^{k-1}-1}{q-1}=q^{k-1}. \qedhere
  \end{equation*}
\end{proof}

In light of Observation~\ref{obs:edges}, determining the numbers of
looped and nonlooped vertices in $G$ is equivalent to finding the
numbers of absolute points of the polarities of $PG(k-1,q)$.

\begin{theorem}\label{thm:num-nonlooped-qeven}
  Let $\Fld_q$ be a finite field having characteristic 2.  One graph in
  $\g_k(\Fld_q)$ will have $\frac{q^{k-1}-1}{q-1}$ nonlooped vertices.  If
  $k$ is even, then the additional graph in $\g_k(\Fld_q)$ will have all
  nonlooped vertices.
\end{theorem}
 
\begin{proof}
  In a field of characteristic~2, since
\begin{align*}
    x^tBx = \sum_{i,j} b_{ij}x_ix_j=\sum_i
    b_{ii}x_i^2+\sum_{i<j}b_{ij}(x_ix_j+x_ix_j) = \sum_i
    b_{ii}x_i^2=\left(\sum_i \sqrt{b_{ii}}x_i\right)^2,
  \end{align*}
  a point $\x$ is absolute if and only if $\sum_i \sqrt{b_{ii}}x_i=0$.

  In a pseudo-polarity, the set of absolute points is the
  hyperplane $\sum_i\sqrt{b_{ii}}x_i=0$.  Since a hyperplane of $PG(k-1,q)$ is a
  projective geometry of projective dimension $k-2$, there are
  $\frac{q^{k-1}-1}{q-1}$ nonlooped vertices in this graph.

  In a null polarity, $b_{ii}=0$ for all $i$.  Therefore every vertex
  is nonlooped (i.e., there are $\frac{q^k-1}{q-1}$ nonlooped
  vertices).  A null polarity occurs when $k$ is even.
\end{proof}

For the odd characteristic case, we will directly apply a standard
result in projective geometry about the number of absolute points in
ordinary polarities.

\begin{theorem}[{\cite[Theorem 22.5.1(b)]{hirschfeld-general-geometries}}]
  Let $q$ be odd.  Then the number of absolute points in a polarity in
  $PG(k-1,q)$ is given by:
  \begin{align*}
    \begin{cases}
      \frac{(q^m-1)(q^{m-1}+1)}{q-1} \text{ or }
      \frac{(q^m+1)(q^{m-1}-1)}{q-1} & \text{if $k=2m$ is even}\\
      \frac{q^{2m}-1}{q-1} & \text{if $k=2m+1$ is odd}
    \end{cases}
  \end{align*}
\end{theorem}

\begin{corollary}\label{cor:num-nonlooped-qodd}
  Let $q$ be odd.  If $k=2m$ is even, then the two graphs in
  $\g_k(\Fld_q)$ will have $\frac{(q^m-1)(q^{m-1}+1)}{q-1}$ and $
  \frac{(q^m+1)(q^{m-1}-1)}{q-1}$ nonlooped vertices, respectively.
  If $k=2m+1$ is odd, then the graph in $\g_k(\Fld_q)$ will have
  $\frac{q^{2m}-1}{q-1}$ nonlooped vertices.
\end{corollary}

We conclude by applying a few standard results for polarities over
$PG(2,q)$ (a projective plane) to give results about $\g_3(\Fld_q)$
and the minimum rank problem.  We note that the polarity graphs of
$PG(2,q)$ for any $q$ are the Erd\H{o}s-R{\'e}nyi graphs from extremal
graph theory (see \cite{erdos-renyi-hungarian}, \cite{erdos-renyi}, or
\cite{brown-erdos-renyi}).  For a survey of interesting properties of
the Erd\H{o}s-R{\'e}nyi graphs and their subgraphs, see
\cite{parsons-projective} or \cite[Chapter 3]{williford-dissert}.

\begin{theorem}\label{thm:white-clique}
  If $G\in\g_3(\Fld_q)$, then the nonlooped vertices in $G$ form a
  clique.
\end{theorem}

\begin{proof}
  Suppose that $u$ and $v$ are distinct nonadjacent nonlooped vertices
  in $G$.  Then $u$ and $v$ are absolute vertices and $u\in\sigma(u)\cap\sigma(v)$
  and $v\in\sigma(u)\cap\sigma(v)$.  This is a contradiction since the intersection
  of any two distinct lines in $PG(2,q)$ is a single point.
\end{proof}

If $G\in\g_3(\Fld_q)$, the formulas in
Theorem~\ref{thm:num-nonlooped-qeven} and
Corollary~\ref{cor:num-nonlooped-qodd} imply that $G$ has $q+1$
nonlooped vertices.  This combined with Corollary~\ref{cor:regular}
and Theorem~\ref{thm:white-clique} gives the following corollary.

\begin{corollary}
  If $G\in\g_3(\Fld_q)$, then each nonlooped vertex is adjacent to
  $q$ nonlooped vertices and $q^2-q$ looped
  vertices.
\end{corollary}

Theorem~\ref{thm:white-clique} also gives us the following theorem.

\begin{theorem}\label{thm:multipartite-weak}
  Let $G=K_{s_1,s_2,\ldots,s_n}$, a simple complete multipartite graph.
  If $q\geq n-1$, then $\mr(\Fld_q,G)\leq 3$.
\end{theorem}

\begin{proof}
  Let $G=K_{s_1,s_2,\ldots,s_n}$.  Then $G$ is a blowup graph of
  $K_n$, where each vertex of $K_n$ is nonlooped.  Since the graph in
  $\g_3(\Fld_q)$ contains a clique of $q+1$ nonlooped vertices, if
  $q+1\geq n$, then $G$ is a blowup graph of the graph in
  $\g_3(\Fld_q)$.
\end{proof}

We can now construct an interesting family of simple graphs.

\begin{theorem}
  For every integer $n\geq1$, let $G_n$ be a simple complete multipartite
  graph $H_1\lor H_2\lor \cdots \lor H_n$ where each $H_i$ is an independent set with
  $s_i>(n-1)^2$ vertices.  We then have $\mr(\Fld_q,G_n)\leq3$ if and
  only if $q\geq n-1$.
\end{theorem}

\begin{proof}
  If $q\geq n-1$, then $\mr(\Fld_q,G_n)\leq3$ by
  Theorem~\ref{thm:multipartite-weak}.

  Conversely, let $q<n-1$.  Let $I$ be the graph in $\g_3(\Fld_q)$ and
  let $I_1$ and $I_2$ be the subgraphs of $I$ induced by the looped
  and nonlooped vertices of $I$, respectively.  Since $I_1$ has $q^2$
  vertices, any blowup of $I_1$ containing more than $q^2$ vertices
  will contain an edge by the pigeon-hole principle.  Since the
  vertices in each $H_i$ form an independent set of size
  $s_i>(n-1)^2>q^2$, at least one vertex in each $H_i$ must be a
  blowup of a vertex in $I_2$.  Furthermore, since the vertices of
  each $H_i$ have the same neighbors, we can assume without loss of
  generality that all of the vertices of each $H_i$ are blowups of
  vertices of $I_1$.  Thus $G_n$ is a blowup of $I_2$.  However, any
  blowup of $I_2$ will be of the form $K_{t_1,t_2,\ldots,t_{q+1}}$
  since $I$ has $q+1$ nonlooped vertices, but $G_n$ is not of this
  form since $q+1<n$.
\end{proof}

\section{Conclusion}

We have suceeded in classifying the structure of graphs in
$\G_k(\Fld_q)$ for any $k$ and any $q$.  We have also shown how this
classification relates to projective geometry.  We have applied a few
results of projective geometry to give results in the minimum rank
problem.

We conclude with a short list of open questions and topics for further
investigation.  First, there are many results about polarity graphs
that could potentially yield results for the minimum rank problem.
What other facts from projective geometry can be applied to give
results in the minimum rank problem over finite fields?

The structural characterization in this paper gives rise to a
theoretical procedure for determining the minimum rank of any graph
over a finite field.  How can this procedure be efficiently
implemented? How can the results of Ding and Kotlov
\cite{ding-kotlov-minrank-finite} be combined with the classification
in this paper to yield results on minimal forbidden subgraphs
describing $\G_k(\Fld_q)$?  The author has implemented such an
algorithm and has some preliminary results on the numbers of forbidden
subgraphs describing $\G_k(\Fld_q)$ for different values of $k$ and
$q$.

Finally, there is still ongoing research investigating the structure
of polarity graphs.  For example, Jason Williford
\cite{williford-dissert}, Michael Newman, and Chris Godsil
\cite{godsil-2005-} have recently investigated the sizes of
independent sets in polarity graphs.  Are there results in the minimum
rank problem that would aid in answering questions about the structure
of polarity graphs?

\section{Acknowledgment}

The author thanks Wayne Barrett and Don March for work in the early
part of this research, including the proof of Theorem~\ref{thm:f2r2}
and early computational experiments, as well as Willem Haemers for
pointing out that the graph $F2R3$ is related to the Fano projective
plane, which led to the investigation of the connection to projective
geometry.  Most of this research was done for the author's
Ph.D.~dissertation and the support of Brigham Young University is
gratefully acknowledged.

\bibliographystyle{alpha}
\bibliography{minrank-proj}

\newcommand{\etalchar}[1]{$^{#1}$}
\begin{thebibliography}{BvdHL05}

\bibitem[AHK{\etalchar{+}}05]{arav}
Marina Arav, Frank~J. Hall, Selcuk Koyuncu, Zhongshan Li, and Bhaskara Rao.
\newblock Rational realizations of the minimum rank of a sign pattern matrix.
\newblock {\em Linear Algebra Appl.}, 409:111--125, 2005.

\bibitem[Alb38]{albert-matrices}
A.~Adrian Albert.
\newblock Symmetric and alternate matrices in an arbitrary field. {I}.
\newblock {\em Trans. Amer. Math. Soc.}, 43(3):386--436, 1938.

\bibitem[Ban07]{bank-thesis}
Efrat Bank.
\newblock Symmetric matrices with a given graph.
\newblock Master's thesis, Technion--Israel Institute of Technology, February
  2007.

\bibitem[BD05]{bento-duarte-tridiag-matrices}
Am{\'e}rico Bento and Ant{\'o}nio~Leal Duarte.
\newblock On {F}iedler's characterization of tridiagonal matrices over
  arbitrary fields.
\newblock {\em Linear Algebra Appl.}, 401:467--481, 2005.

\bibitem[BF07]{barioli-fallat-joins}
Francesco Barioli and Shaun Fallat.
\newblock On the minimum rank of the join of graphs and decomposable graphs.
\newblock {\em Linear Algebra Appl.}, 421(2-3):252--263, 2007.

\bibitem[BFH04]{BFH1-minrankpath}
Francesco Barioli, Shaun Fallat, and Leslie Hogben.
\newblock Computation of minimal rank and path cover number for certain graphs.
\newblock {\em Linear Algebra Appl.}, 392:289--303, 2004.

\bibitem[BFH05a]{bfh-mult-pathcover-tree}
Francesco Barioli, Shaun Fallat, and Leslie Hogben.
\newblock On the difference between the maximum multiplicity and path cover
  number for tree-like graphs.
\newblock {\em Linear Algebra Appl.}, 409:13--31, 2005.

\bibitem[BFH05b]{bfh-cdv}
Francesco Barioli, Shaun Fallat, and Leslie Hogben.
\newblock A variant on the graph parameters of {C}olin de {V}erdi\`ere:
  implications to the minimum rank of graphs.
\newblock {\em Electron. J. Linear Algebra}, 13:387--404 (electronic), 2005.

\bibitem[BGL]{barrett-grout-loewy-mrF2R3}
Wayne Barrett, Jason Grout, and Raphael Loewy.
\newblock The minimum rank problem over the finite field of order 2: minimum
  rank 3.
\newblock 38 pages. Preprint available at
  \mbox{\url{http://arxiv.org/abs/math.CO/0612331}}.

\bibitem[Bro66]{brown-erdos-renyi}
W.~G. Brown.
\newblock On graphs that do not contain a {T}homsen graph.
\newblock {\em Canad. Math. Bull.}, 9:281--285, 1966.

\bibitem[BvdHL04]{barrett-vdHL-minrank2-infinite}
Wayne Barrett, Hein van~der Holst, and Raphael Loewy.
\newblock Graphs whose minimal rank is two.
\newblock {\em Electron. J. Linear Algebra}, 11:258--280 (electronic), 2004.

\bibitem[BvdHL05]{barrett-vdHL-minrank2-finite}
Wayne Barrett, Hein van~der Holst, and Raphael Loewy.
\newblock Graphs whose minimal rank is two: the finite fields case.
\newblock {\em Electron. J. Linear Algebra}, 14:32--42 (electronic), 2005.

\bibitem[CDH{\etalchar{+}}07]{hogben-minrank-tree}
Nathan~L. Chenette, Sean~V. Droms, Leslie Hogben, Rana Mikkelson, and Olga
  Pryporova.
\newblock Minimum rank of a tree over an arbitrary field.
\newblock {\em Electron. J. Linear Algebra}, 16:183--186 (electronic), 2007.

\bibitem[CdV98]{cdv-eig-mult}
Yves Colin~de Verdi{\`e}re.
\newblock Multiplicities of eigenvalues and tree-width of graphs.
\newblock {\em J. Combin. Theory Ser. B}, 74(2):121--146, 1998.

\bibitem[CG01]{gr-graph-theory}
Gordon~Royle Chris~Godsil.
\newblock {\em Algebraic Graph Theory}.
\newblock Springer-Verlag, 2001.

\bibitem[CHLW03]{chen}
Guantao Chen, Frank~J. Hall, Zhongshan Li, and Bing Wei.
\newblock On ranks of matrices associated with trees.
\newblock {\em Graphs Combin.}, 19(3):323--334, 2003.

\bibitem[Coh03]{cohn-algebra}
P.~M. Cohn.
\newblock {\em Basic algebra}.
\newblock Springer-Verlag London Ltd., London, 2003.
\newblock Groups, rings and fields.

\bibitem[DK06]{ding-kotlov-minrank-finite}
Guoli Ding and Andre{\u\i} Kotlov.
\newblock On minimal rank over finite fields.
\newblock {\em Electron. J. Linear Algebra}, 15:210--214 (electronic), 2006.

\bibitem[ER62]{erdos-renyi-hungarian}
P.~Erd{\H{o}}s and A.~R{\'e}nyi.
\newblock On a problem in the theory of graphs.
\newblock {\em Magyar Tud. Akad. Mat. Kutat\'o Int. K\"ozl.}, 7:623--641
  (1963), 1962.

\bibitem[ERS66]{erdos-renyi}
P.~Erd{\H{o}}s, A.~R{\'e}nyi, and V.~T. S{\'o}s.
\newblock On a problem of graph theory.
\newblock {\em Studia Sci. Math. Hungar.}, 1:215--235, 1966.

\bibitem[GN05]{godsil-2005-}
C~D Godsil and M~W Newman.
\newblock Eigenvalue bounds for independent sets, 2005.

\bibitem[Hir98]{hirschfeld-projective}
J.~W.~P. Hirschfeld.
\newblock {\em Projective geometries over finite fields}.
\newblock Oxford Mathematical Monographs. The Clarendon Press Oxford University
  Press, New York, second edition, 1998.

\bibitem[HLR04]{hall}
Frank~J. Hall, Zhongshan Li, and Bhaskara Rao.
\newblock From {B}oolean to sign pattern matrices.
\newblock {\em Linear Algebra Appl.}, 393:233--251, 2004.

\bibitem[Hsi01]{hsieh-minrank}
Liang-Yu Hsieh.
\newblock {\em On Minimum Rank Matrices having a Prescribed Graph}.
\newblock PhD thesis, University of Wisconsin, Madison, 2001.

\bibitem[HT91]{hirschfeld-general-geometries}
J.~W.~P. Hirschfeld and J.~A. Thas.
\newblock {\em General {G}alois geometries}.
\newblock Oxford Mathematical Monographs. The Clarendon Press Oxford University
  Press, New York, 1991.
\newblock , Oxford Science Publications.

\bibitem[JD99]{johnson-duarte-trees}
C.~R. Johnson and Ant{\'o}nio~Leal Duarte.
\newblock The maximum multiplicity of an eigenvalue in a matrix whose graph is
  a tree.
\newblock {\em Linear and Multilinear Algebra}, 46(1-2):139--144, 1999.
\newblock Invariant factors.

\bibitem[JS02]{johnson-saiago-maxmult}
Charles~R. Johnson and Carlos~M. Saiago.
\newblock Estimation of the maximum multiplicity of an eigenvalue in terms of
  the vertex degrees of the graph of a matrix.
\newblock {\em Electron. J. Linear Algebra}, 9:27--31 (electronic), 2002.

\bibitem[KSS97]{blowup-lemma}
J{\'a}nos Koml{\'o}s, G{\'a}bor~N. S{\'a}rk{\"o}zy, and Endre Szemer{\'e}di.
\newblock Blow-up lemma.
\newblock {\em Combinatorica}, 17(1):109--123, 1997.

\bibitem[Nyl96]{nylen-minrank}
Peter~M. Nylen.
\newblock Minimum-rank matrices with prescribed graph.
\newblock {\em Linear Algebra Appl.}, 248:303--316, 1996.

\bibitem[Par76]{parsons-projective}
T.~D. Parsons.
\newblock Graphs from projective planes.
\newblock {\em Aequationes Math.}, 14(1-2):167--189, 1976.

\bibitem[Sin06]{sinkovic-thesis}
John Sinkovic.
\newblock The relationship between the minimal rank of a tree and the
  rank-spreads of the vertices and edges.
\newblock Master's thesis, Brigham Young University, December 2006.

\bibitem[Ste07]{sage-2.8.15}
William Stein.
\newblock {\em {Sage} {M}athematics {S}oftware ({V}ersion 2.8.15)}.
\newblock The Sage~Group, 2007.
\newblock \url{http://www.sagemath.org}.

\bibitem[vdH03]{vdh-nullity}
Hein van~der Holst.
\newblock Graphs whose positive semi-definite matrices have nullity at most
  two.
\newblock {\em Linear Algebra Appl.}, 375:1--11, 2003.

\bibitem[Wil04]{williford-dissert}
Jason Williford.
\newblock {\em Constructions in finite geometry with applications to graphs}.
\newblock PhD thesis, University of Delaware, 2004.

\end{thebibliography}

\appendix

\section{Sage code to generate graphs}
\label{app:gen-subs-graphs}
\begin{lstlisting}
# This code is written for Sage 2.8.15.
# See http://www.sagemath.org/

def bilinearforms(F,mr):
    # Construct a matrix space for our bilinear forms
    MSpace = MatrixSpace(F,mr)
    # The identity matrix is always 
    # a congruence class representative
    forms = [MSpace.identity_matrix()]
    # Add the extra matrices in the even rank cases
    if(mod(mr,2)==0): # even rank
	if(F.characteristic()==2): # characteristic 2
	    # Add $\diag(H_1,H_2,\ldots,H_{\text{mr}/2})$
	    hyperbolic = matrix(F,[[0,1],[1,0]])
            hyperbolic_form = matrix(F,[])
            for i in (1..Integer(mr/2)):
                hyperbolic_form = hyperbolic_form.block_sum(hyperbolic)
            forms.append(hyperbolic_form)
	else: # odd characteristic
	    # Add $\diag(I_{n-1},\nu)$, where $\nu$ is a non-square
	    nonidentity_form = MSpace.identity_matrix()
	    # Find a non-square
	    for nu in F:
		if not nu.is_square():
		    break
	    if nu.is_square():
                raise NotImplementedError, \
                      "Cannot find a non-square in the field."
	    nonidentity_form[mr-1, mr-1] = nu
	    forms.append(nonidentity_form)
    return forms

def get_matrices(F,mr):
    # $U$ has one vector for every equivalence class in $PG(mr-1,q)$
    U = matrix(F,[list(v) 
		for v in ProjectiveSpace(mr-1,F)]).transpose()
    B = bilinearforms(F,mr)
    return U, B

def get_graphs(F,mr):
    U, B = get_matrices(F,mr)
    product_matrices = [U.transpose()*b*U for b in B]
    graphs = [Graph(m) for m in product_matrices]
    for i in range(len(graphs)):
        graphs[i].loops(true);
        graphs[i].add_edges([[j,j] for j in range(len(graphs[i])) \
			     if product_matrices[i][j,j] != 0])
    return graphs

def show_graphs(F,mr):
    for g in get_graphs(F,mr):
	# Vertices with loops are black, others are white
        vcolors={'black': g.loop_vertices(),\
                 'white': [i for i in g.vertices()
                           if i not in g.loop_vertices()]}
        g.show(layout='circular', vertex_colors=vcolors, \
	       vertex_labels=false)

# To retrieve the matrices for graphs in $\g_3(\Fld_2)$:
U, B = get_matrices(F=FiniteField(2), mr=3)

# To retrieve the graphs in $\g_3(\Fld_2)$:
graph_list = get_graphs(F=FiniteField(2), mr=3)

# To display the graphs with vertices colored appropriately:
show_graphs(F=FiniteField(2), mr=3)
\end{lstlisting} 

\end{document}